\input harvmac
\overfullrule=0pt
%

\def\bar{\overline}
\def\x{{\times}}
\def\ra{{\rightarrow}}

\Title{ \vbox{\baselineskip12pt
\hbox{math.AG/0012196}
\hbox{HUB-EP-00/60}}}
{\vbox{\centerline{Fourier-Mukai Transform and Mirror Symmetry}
\bigskip\centerline{for D-Branes on Elliptic Calabi-Yau}}}
\centerline{Bj\"orn Andreas$^{1}$, Gottfried Curio$^{2}$,
Daniel Hern\'andez Ruip\'erez$^{3}$ and Shing-Tung Yau$^{1}$}
\bigskip
\centerline{\it $^{1}$Department of Mathematics,
Harvard University, Cambridge, MA 02138, USA}
\smallskip
\centerline{\it $^2$Humboldt-Universit\"at zu Berlin,
Institut f\"ur Physik, D-10115 Berlin, Germany}
\smallskip
\centerline{\it $^{3}$Departamento de Matem\'aticas,
Universidad de Salamanca,  37008, Salamanca, Spain}

\baselineskip 18pt
\def\sqr#1#2{{\vbox{\hrule height.#2pt\hbox{\vrule width
.#2pt height#1pt \kern#1pt\vrule width.#2pt}\hrule height.#2pt}}}

\smallskip
\smallskip
\smallskip
\noindent
Fibrewise T-duality (Fourier-Mukai transform)
for D-branes on an elliptic Calabi-Yau three-fold
$X$ is seen to have an expected adiabatic form for its induced
cohomology operation only when an appropriately twisted operation
resp. twisted charge is defined. Some differences
with the case of $K3$ as well as connections with the spectral cover
construction for bundles on $X$ are pointed out. In the context of
mirror symmetry Kontsevich's association of line bundle twists
(resp. a certain 'diagonal' operation) with monodromies (esp. the
conifold monodromy) is made explicit and checked for two example models.
Interpreting this association as a relation between FM transforms
and monodromies, we express
the fibrewise FM transform through known
monodromies. The operation of fibrewise duality as well as the
notion of a certain index relevant to the computation of the moduli
space of the bundle is transported to the sLag side. Finally
the moduli space for D4-branes and its behaviour under the FM
transform is considered with an application to the spectral cover.

\Date{}

\lref\hart{R. Hartshorne, ``Algebraic Geometry'', Springer 1977.}

\lref\RD{R. Hartshorne, ``Residues and Duality'', Springer 1966}

\lref\Karo{F. Hirzebruch, ``Topological Methods in Algebraic Geometry'',
Springer 1978\semi
M. Karoubi, ``K-Theory'', Springer 1978.}

\lref\wik{E. Witten, ``D-Branes and K-theory, hep-th/9810188.}

\lref\koba{S. Kobayashi, ``Differential Geometry of Complex Vector Bundles'',
Princeton University Press 1987.}

\lref\donal{S. K. Donaldson, ``A polynomial Invariant for Smooth
Four-Manifolds'', Topology {\bf 29} (1990) 257.}

\lref\thom{R. P. Thomas, ``A Holomorphic Casson Invariant For Calabi-Yau
Three-Folds, and Bundles on K3 Fibrations'', math.AG/9806111.}

\lref\FMW{R. Friedman, J. Morgan and E. Witten, ``Vector Bundles and F-
Theory,'' Commun. Math. Phys. {\bf 187} (1997) 679, hep-th/9701162.}

\lref\cur{G. Curio, ``Chiral matter and transitions in heterotic string
models'', Phys. Lett. {\bf B435} (1998) 39, hep-th/9803224.}

\lref\CDo{G. Curio and R. Donagi, ``Moduli in {N=1} Heterotic/F-Theory
 Duality'', Nucl.Phys. {\bf B518} (1998) 603, hep-th/9801057.}

\lref\vertfivebr{B. Andreas and G. Curio, ``Horizontal and Vertical
Five-Branes in Heterotic/F-Theory Duality'', JHEP 0001 (2000) 013,
hep-th/9912025.}

\lref\Euler{B. Andreas and G. Curio, ``On discrete Twist
and Four-Flux in N=1 heterotic/F-theory compactifications'',
hep-th/9908193.}

\lref\ack{B. Andreas, G. Curio and A. Klemm, ``Towards the Standard
Model spectrum from elliptic Calabi-Yau'',
hep-th/9903052.}

\lref\fivebr{B. Andreas and G. Curio, ``Three-Branes
and Five-Branes in N=1 Dual String Pairs'', Phys.Lett. {\bf B417} (1998) 41,
hep-th/9706093.}

\lref\chir{G. Curio, ``Chiral Multiplets
in N=1 Dual String Pairs'', Phys.Lett. {\bf B409} (1997) 185,
hep-th/9705197.}

\lref\acl{B. Andreas, G. Curio and D. Lust, ``N=1 Dual String Pairs
and their Massless Spectra'', Nucl.Phys. {\bf B507} (1997) 175,
hep-th/9705174.}

\lref\bax{B. Andreas, ``On Vector Bundles and Chiral Matter in N=1
Heterotic Compactifications'', JHEP {\bf 9901} (1999) 011, hep-th/9802202}

\lref\micha{M. Bershadsky, A. Johansen, T. Pantevand V. Sadov,``On
Four-Dimensional Compactifications of F-Theory'', Nucl.Phys.
{\bf B505} (1997) 165, hep-th/9701165.}

\lref\diacrom{D.-E. Diaconescu and C. R\"omelsberger, ``D-Branes and Bundles
on Elliptic Fibrations'', hep-th/9910172.}

\lref\cand{P. Candelas, X.C. de la Ossa, P.S. Green and L. Parkes,
``A pair of Calabi-Yau manifolds as an exactly solvable superconformal
theory'', Nucl. Phys. {\bf B359} (1991) 21.}

\lref\morri{D.R. Morrison, ``Geometric Aspects of Mirror Symmetry'',
math.AG/0007090.}

\lref\kllw{P. Kaste, W. Lerche, C. A. L\"utken and J. Walcher,
``D-Branes on K3-Fibrations'', hep-th/9912147.}

\lref\morr{P. Candelas, A. Font, S. Katz and D. Morrison, ``Mirror Symmetry
for Two Parameter Models- II'', hep-th/9403187.}

\lref\acm{B. Andreas, G. Curio and R. Minasian, ``Anomalous D-Brane Charge in
F-Theory Compactifications'', hep-th/0007212, JHEP 0009 (2000) 022}

\lref\aty{M. F. Atiyah and R. Bott, ``A Lefschetz Fixed Point Formula for
Elliptic Differential Operators'', Bull. Am. Math. Soc. {\bf 72} (1966) 245.}

\lref\att{M. F. Atiyah and G. B. Segal, ``The Index of Elliptic Operators,
II'', Ann. Math. {\bf 87} (1968) 531.}

\lref\bridge{T. Bridgeland, `` Fourier-Mukai transforms for elliptic
surfaces'', J. reine angew. Math. {\bf 498}, alg-geom/9705002.}

\lref\dave{P. Candelas, X. De La Ossa, A. Font, S. Katz and D. R. Morrison,
``Mirror symmetry
for two parameter models I'', Nucl. Phys. {\bf B416} (1994) 481,
hep-th/9308083.}

\lref\hoskl{S. Hosono, A. Klemm, S. Theissen and S. T. Yau, ``Mirror
symmetry, mirror map and
applications to complete intersection Calabi-Yau spaces'',
Nucl. Phys. {\bf B433} (1995) 501, hep-th/9406055}

\lref\grepl{B. R. Greene and M. R. Plesser, ``Duality in Calabi-Yau moduli
spaces'', Nucl. Phys. {\bf B338} (1990) 15}

\lref\bart{C. Bartocci, U. Bruzzo, D. Hern\'andez Ruip\'erez and J.
M. Mu\~noz  Porras, ``Mirror symmetry on K3 surfaces via
Fourier-Mukai transform'', Commun. Math. Phys. {\bf 195} (1998), 79-93}

\lref\RPo{D. Hern\'andez Ruip\'erez and J. M. Mu\~noz Porras,
``Structure of the  moduli space of stable sheaves on elliptic
fibrations'', Preprint Universidad de Salamanca, (1998)}

\lref\gm{S. I. Gelfand and Y. I. Manin, ``Methods of Homological
Algebra'',  Springer 1996}

\lref\Thomas{R. P. Thomas, ``Mirror symmetry
and actions of braid groups on derived categories'',
math.AG/0001044.}

\lref\hosono{S. Hosono, ``Local Mirror Symmetry and
Type IIA Monodromy of Calabi-Yau manifolds'',
Adv. Theor. Math. Phys. {\bf 4} (2000), hep-th/0007071.}

\lref\zaslow{A. Polishchuk and E. Zaslow,
``Categorical Mirror Symmetry: The Elliptic Curve'',
Adv. Theor. Math. Phys. {\bf 2} (1998) 443, math.AG/9801119.}

\lref\AspDon{P.S. Aspinwall and R.Y. Donagi, ``The Heterotic String,
the Tangent Bundle, and Derived Categories'',
Adv. Theor. Math. Phys. {\bf 2} (1998) 1041,
hep-th/9806094.}

\lref\Muk{S. Mukai, ``Duality between $D(X)$ and $D(\hat X)$ with
its application to Picard sheaves'',  Nagoya Math.~J.~{\bf 81} (1981),
153.}

\lref\Mac{A. Maciocia, ``Generalized Fourier-Mukai Transforms'',
J. reine angew. Math., {\bf 480} (1996) 197,
alg-geom/9705001.}

\lref\BridgeMac{T. Bridgeland and A. Maciocia, ``Fourier-Mukai
transforms for K3 and elliptic fibrations'',
math.AG/9908022.}

\lref\JarMac{M. Jardim and A. Maciocia, ``A Fourier-Mukai approach
to spectral data for instantons'',
math.AG/0006054.}

\lref\LeuYauZas{N.C. Leung, S.-T. Yau and E. Zaslow, ``From
Special Lagrangian to Hermitian-Yang-Mills via Fourier-Mukai Transform'',
math.DG/0005118.}

\lref\SeiTho{P. Seidel and R.P. Thomas, ``Braid group actions
on derived categories of coherent sheaves'',
math.AG/0001043.}

\lref\Man{Y.I. Manin, ``Moduli, Motives, Mirrors'',
math.AG/0005144.}

\lref\mdoug{M.R. Douglas, ``Topics in D-Geometry'', hep-th/9910170.}

\lref\doug{M.R. Douglas, ``D-Branes on Calabi-Yau Manifolds'',
math.AG/0009209.}

\lref\DouFioRoem{M.R. Douglas, B. Fiol and C. R\"omelsberger,
``Stability and BPS branes'',
hep-th/0002037.}

\lref\DouFioRoemsec{M.R. Douglas, B. Fiol and C. R\"omelsberger,
``The spectrum of BPS branes on a noncompact Calabi-Yau'',
hep-th/0003263.}

\lref\Dou{M.R. Douglas, ``D-branes, Categories and N=1
Supersymmetry'', hep-th/0011017.}

\lref\hor{R.P. Horja, ``Hypergeometric functions and
mirror symmetry in toric varieties'',
math.AG/9912109.}

\lref\DiaGom{D.E. Diaconescu and J. Gomis, ``Fractional Branes
and Boundary States in Orbifold Theories'', JHEP 0010 (2000) 001,
hep-th/9906242.}

\lref\FioMar{B. Fiol and M. Marino, ``BPS states and algebras from quivers'',
JHEP 0007 (2000) 031,
hep-th/0006189.}

\lref\KachMcGre{S. Kachru and J. McGreevy, ``Supersymmetric
Three-cycles and (Super)symmetry Breaking'', Phys.Rev. D61 (2000)
026001, hep-th/9908135.}

\lref\DOPWa{R. Donagi, B. Ovrut, T. Pantev and D. Waldram,
``Spectral involutions on rational elliptic surfaces'',
math.AG/0008011.}

\lref\DOPWb{R. Donagi, B. Ovrut, T. Pantev and D. Waldram,
``Standard-model bundles'',
math.AG/0008010.}

\lref\DOPWc{R. Donagi, B. Ovrut, T. Pantev and D. Waldram,
``Standard-Model Bundles on Non-Simply Connected Calabi--Yau Threefolds'',
hep-th/0008008.}

\lref\GraGurOv{A. Grassi, Z. Guralnik and B. Ovrut,
``Five-Brane BPS States in Heterotic M-Theory'',
hep-th/0005121.}

\lref\DOPWd{R. Donagi, B. Ovrut, T. Pantev and D. Waldram,
``Standard Model Vacua in Heterotic M-Theory'',
hep-th/0001101.}

\lref\DOPWe{R. Donagi, B. Ovrut, T. Pantev and D. Waldram,
``Standard Models from Heterotic M-theory'',
hep-th/9912208.}

\lref\LukOv{A. Lukas and B.A. Ovrut, ``Symmetric Vacua in Heterotic M-Theory'',
hep-th/9908100.}

\lref\DOW{R. Donagi, B. Ovrut and D. Waldram,
``Moduli Spaces of Fivebranes on Elliptic Calabi-Yau Threefolds'',
JHEP 9911 (1999) 030, hep-th/9904054.}

\lref\DLOW{R. Donagi, A. Lukas, B. Ovrut and D. Waldram,
``Non-Perturbative Vacua and Particle Physics in M-Theory'',
JHEP 9905 (1999) 018, hep-th/9811168.}

\lref\SYZ{A. Strominger, S.T. Yau and E. Zaslow,
``Mirror Symmetry is T-Duality'', Nucl.Phys. {\bf B479} (1996) 243,
hep-th/9606040.}

\lref\Yo{K. Yoshioka, ``Moduli spaces of stable sheaves on
abelian surfaces'',  math.AG/0009001.}


\newsec{Introduction and summary}

This paper treats some connections between four different, although
related topics: D-branes, mirror symmetry, elliptic Calabi-Yau and
Fourier-Mukai transform.

The last year saw an intense study on BPS {\it D-branes} in type II string
theories on a Calabi-Yau manifold
(cf. \wik ,\mdoug ,\doug ,\DouFioRoem ,\DouFioRoemsec
,\Dou ,\diacrom ,\DiaGom ,\kllw ,\FioMar ,\KachMcGre ),
focussing on the behaviour of the D-brane spectrum under variation in
the Calabi-Yau moduli space (related to the phenomenon of marginal
stability) as well as on the relations at a special point in moduli
space (such as the relation with boundary conformal field theory).

This development has some close connections with a reformulation of
{\it mirror symmetry} given by Kontsevich on the one hand and
Strominger/Yau/Zaslow\foot{Note that the T-duality on the $T^3$ fibre
in that construction is not easily related to the T-duality on
the holomorphic elliptic fibre considered later in the framework of
the fibrewise FM transform as the holomorphic $T^2$ is not contained
in the $T^3$.} and Vafa
on the other which brings supersymmetric D-branes
on both sides of the mirror correspondence into the play
(cf. \cand ,\SYZ ,\morri ,\dave ,\morr ,\hoskl ,\grepl ,\hosono
,\zaslow ,\LeuYauZas ,\SeiTho ,\Man ,\hor ).
This relates bundles or
better sheaves (think of these here as bundles supported on
holomorphic subvarieties) on a Calabi-Yau $X$ in type IIA string
theory with special Lagrangian submanifolds (with an $U(1)$ bundle over
them) in the mirror Calabi-Yau $Y$ (or more precisely the derived
category $D(X)$ of the category of sheaves on $X$ with Fukaya's
$A_{\infty}$ category of Lagrangian submanifolds of $Y$);
cohomological invariants relate
then the $H^{even}(X)$ of the bundle side and $H^3(Y)$ on the sLag
side. Furthermore monodromies around divisors ${\cal D}$ in moduli space where
some such even-dim. cycle (say a divisor $D$) vanish correspond
conjecturally to twisting with the line bundle ${\cal L}_D$
associated to the divisor $D$ and a similar more complicated relation,
which generalizes the twisting with the help of the concept of a
Fourier-Mukai transform, relates the conifold monodromy.

Now for the description of bundles on a Calabi-Yau it was taken a
great step forward when Friedman/Morgan/Witten made explicit (via spectral
covers) a construction of bundles for the case of an
{\it elliptically fibered Calabi-Yau}. This was then intended as
compactification space for the heterotic string and allowed detailed
studies on such issues as the moduli space, brane impurities, relation
with F-theory and model building
(cf. \FMW ,\micha ,\cur ,\bax ,\acl ,\chir ,\fivebr ,\CDo ,\acm ,\ack
,\Euler ,\vertfivebr ,\DLOW ,\DOW ,\LukOv ,\DOPWe ,\DOPWd ,\DOPWc
,\DOPWb ,\DOPWa ,\GraGurOv ).

This class of elliptic Calabi-Yau's is also especially interesting as
it allows for a version of the {\it Fourier-Mukai transform} of a bundle $V$
on non-toroidal spaces which in contrast to earlier transforms in such
cases (on $K3$ say) keeps completely the idea of using the duality on
a torus by building a fibrewise FM transform
(cf. \bridge ,\bart ,\RPo ,\Thomas ,\AspDon ,\Mac ,\BridgeMac ,\JarMac ).
This is what the
physicist would call T-duality on the fibre and operates in an
interesting way on the spectrum of D-branes on the one hand and has on
the other hand a close connection with the spectral cover construction
(a relation of an FM transform to Kontsevich's version of mirror
symmetry was already mentioned).

A number of related points in the aforementioned web of connections
will be studied in the paper.
In {\it section 2} the cohomological invariants of a bundle and its
fibrewise dual are computed\foot{with (appendix) and without (main
body of the paper) use of the spectral cover construction} and it is
shown how by using an appropriately twisted charge (in analogy with
the $\sqrt{Td(X)}$ twist of the Chern character to get the Mukai
vector in the usual full FM transform) the adiabatic
character of the operation is confirmed. I.e. by using a decomposition
of the cohomology into base and fibre parts the operation of the
fibrewise duality on the cohomology will be seen to take the form
one gets from an adiabatic extension of the same operation on the
cohomology of a $T^2$ (fulfilling the expectations from the
interpretation as T-duality on D-branes); a simpler variant of the
twisting idea is seen also to be necessary for the case of $K3$ where,
in accordance with earlier treatments of that case in the literature,
still an operation relating the untwisted Chern characters themselves
can be given by considering a natural twist in the duality functor
itself; after these cases of Calabi-Yau's of complex dimension one and
two our presentation follows the line of ascending complexity and
demonstrates how for a Calabi-Yau three-fold one has to use the
twisted charge
definition (which also naturally incorporates via an reinterpretation
the findings for $K3$) or an adapted version of Mukai's {\bf f}-map
(which uses the usual Mukai vector as charge but a slightly twisted operation).

In {\it section 3} we treat some well-studied two parameter CY
(represented by hypersurfaces of degree 8 and 12 in weighted
projective space) and make explicit Kontsevich's association of
monodromies with twists by line bundles resp. a more complicated
operation for the conifold monodromy.

In {\it section 4} a connection between the two main themes of the
foregoing chapters is given: the fibrewise FM transform is given its
place in Kontsevich's general association of Fm transforms with
monodromies. For this note that
the relative FM transform on elliptic fibrations is an
autoequivalence $S\colon D(X)\to D(X)$ of the derived category,
its inverse\foot{To be precise, one should use here the term
quasi-inverse instead of inverse  and the equality signs in
$\hat{S}\circ S (\cdot)=(\cdot)[-1]$ and
$S\circ\hat{S}(\cdot)=(\cdot)[-1]$ have to be understood as natural
functor isomorphisms (see pag. 71 of \gm\ ; there is also given
a quite illuminating discussion about
how inappropriate the notion of isomorphism of categories is).}
functor being (up to a shift) the functor $\hat{S}$
described in Section 2.2.
Then $S$ should correspond in the mirror to a monodromy.
As part of Kontsevich's generalization of mirror symmetry can be
formulated without the mirror\foot{the twists with the line bundle
associated to a divisor are related to
the monodromies around the locus in moduli space where the divisor
vanishes, but this already in the K\"ahler moduli space of type IIA,
i.e. $H^{even}(X)$  -  the identification of the period monodromy in
the complex structure moduli space of the mirror CY $Y$ (related to $H^3(Y)$)
are then reached by combining the first identification with the mirror
map which identifies the monodromies} we can test this already on $X$.

In {\it section 5} we point to another connection between our fibrewise
FM transform and mirror symmetry. Namely it is pointed out how via the
mirror identification one can transport operation of fibrewise FM on
cohomology studied in section 2 to the mirror side and get there a
corresponding operation on the middle cohomology (this is basically
just a careful comparison of base choices as one has to relate natural
symplectic bases of periods with the decomposition into base and fibre
parts used to make the fibrewise Fm transform most easily visible). It
is then of course a very interesting question whether that operation
on $H^3(Y)$ transported from $H^{even}(X)$ actually can be derived
from a certain operation {\it already on the space level} as it can be done
on the original bundle side. As a second instance of that transport
philosophy we point out how the more precise information one has on
the moduli space of bundles in the elliptic set-up (as given
essentially by a certain index) can be transported to a "sLag index"
on the sLag side which should give more structure to conjectured
relations such as the one betwen $h^1(End V)$ and $h^1(Q)$ ($Q$ the
sLag 3-cycle).

Finally in {\it section 6} again some moduli space questions in connection with
$h^1(End V)$ on the bundle side
are treated; this shows how these
quantities are related to corresponding expressions in the spectral
cover construction; as that construction relates D6-brane (the bundle
$V$) to a D4-brane (its spectral cover divisor
with a certain line bundle over it) this is considered from
the general perspective of studying the action of the fibrewise FM
transform relating D6-branes to D4-branes on $X$.

The appendix recalls the computation of cohomological invariants of a
bundle $V$ and its fibrewise dual in the spectral cover construction
and shows why the twist of the duality functor alone, which was
successful in the $K3$ case to get the 'adiabatic' transformation
matrix, is insufficient in the three-fold case.

Let us state the technical framework the paper will be moving in.
We will assume the
elliptically fibered CY $X$ has a smooth Weierstrass model, having
singular $I_1$ fibers over a one-dimensional locus in the base, which
furthermore has a section $\sigma$. All sheaves are coherent. Quite often
we will
deal with bundles which are fibrewise of degree zero. They will be in
addition semistable on the generic fibre in case we refer to the spectral cover
construction. For those sheaves the fibrewise FM transforms preserves fibrewise
semistability. The preservation under fibrewise FM of the absolute
stability (with
respect to some particular kind of polarization in $X$) is only known in
the two
dimensional case (cf. \RPo ,\JarMac ,\Yo ). The natural normalization of the
Poincar\'e bundle we use will be recalled in section 2. As usual in
consideration of the FM transform the 'dual' CY
${\tilde X}$, given as a compactification of the Jacobian of the original
fibration
(Jacobian fibrewise), will be identified with $X$ when appropriate.
Concerning the
notion of a Dp-brane which has its p-dimensional spatial world-volume
wrapped on a
holomorphic p-cycle in $X$ the mathematical oriented reader should
think of a bundle concentrated on that p-cycle, i.e. a sheaf on $X$
with support on this cycle (cf. subsection 2.5
and the presentations \Dou ,\doug ,\mdoug).

\newsec{Fibrewise Fourier-Mukai transform (T-duality) on elliptic Calabi-Yau}
We consider an $SU(n)$ bundle $V$ over an elliptically fibered
$X$ of $c_1(V)=0$ or equivalently
$n$ D6-branes wrapped over $X$ with induced lower-dimensional
D2-and D0-brane charges
(D2i-charges meaning here for now just $ch_i(V)$).
The case of elliptically fibered Calabi-Yau three-folds $X$ will have a
double advantage: one can describe $V$ explicitly
via the associated spectral cover \FMW\ and has
furthermore the action of T-duality on the elliptic fibre on the
bundles.
So our procedure in this section will be first to make the
bundle description more explicit using the spectral cover method and
then to describe the change in $ch(V)$ induced by fibrewise T-duality, an
operation which should be mirrored on the sLag side by a corresponding
operation.

The T-duality on the $T^2$ fiber maps in general
(the subscripts indicate whether fibre $F$ or base $B$ is contained
(resp. contains) the wrapped world-volume)
\eqn\Dtual{\eqalign{ D6&\rightarrow \tilde{D4}_B\cr
D4_B\rightarrow \tilde{D6}\;\; &,
\;\; D4_F\rightarrow\tilde{D2}_B \cr
   D2_B\rightarrow \tilde{D4}_{\tilde{F}}\;\; &,
\;\; D2_F\rightarrow  \tilde{D0}\cr
                    D0&\rightarrow \tilde{D2}_{\tilde{F}}}}
Our goal in this section is to understand the operation of fibrewise T-duality
on the cohomological data representing the bundle $V$ and its
Fourier-Mukai (FM) dual $\tilde{V}$, i.e. we will mirror the mentioned
D-brane relations as a map between $ch(V)$ and
$ch(\tilde{V})$ (in a first approximation; the modification needed to
make this work will be made precise along the way). For this
we will assume the following decomposition of the vertical cohomology
(the ${\bf C}$ in $H^4(X)$ resp. $H^6(X)$ are $H^4(B)$ resp. $\sigma
H^4(B)$; the latter $\sigma$ will be often suppressed)
\eqn\decomp{\eqalign{
H^0(X)&={\bf C}\cr
\oplus&\cr
H^2(X)&={\bf C}\sigma \oplus H^2(B)\cr
\oplus&\cr
H^4(X)&=H^2(B)\sigma \oplus {\bf C}\cr
\oplus&\cr
H^6(X)&={\bf C}}}
Then essentially the six entries of the Chern character
vector in our decomposition are
pairwise interchanged and the transformation has the block-diagonal form
\eqn\genmatrix{
Q=\pmatrix{
0&*&&&&\cr
{*}&0&&&&\cr
&&0&*&&\cr
&&*&0&&\cr
&&&&0&*\cr
&&&&*&0} \tilde{Q}}
This has the following interpretation: the fibrewise Fourier-Mukai
transform (fibrewise T-duality) is given here just by adiabatic
extension of the T-duality on the
fibre; this adiabatic relation should be reflected by a corresponding
adiabatic relation between the matrices representing the operation on
the cohomology. Now
for the case of an one-dimensional Calabi-Yau
consisting just of an elliptic curve (representing the fibre $T^2$)
the matrix is \bridge\
\eqn\ellmatrix{
A=\pmatrix{
0&1\cr
-1&0}}
So one would like to see that
the actual form of the transformation matrix
induced by the mathematical operation of fibrewise Fourier-Mukai transformation
on the bundle (which represents the fibrewise T-duality process)
is given by
(in the order of
arrangement given in the decomposition above)
\eqn\precisematrix{
Q=\pmatrix{0&1&0&0&0&0\cr
-1&0&0&0&0&0\cr
0&0&0&1&0&0\cr
0&0&-1&0&0&0\cr
0&0&0&0&0&1\cr
0&0&0&0&-1&0}\tilde{Q}}
thereby then confirming the 'adiabatic' interpretation of the Fourier-Mukai
transform as
\eqn\adiabmatrix{
\pmatrix{
{\bf 0}_3&{\bf 1}_3\cr
-{\bf 1}_3&{\bf 0}_3}}
with entries now consisting of
$3\x 3$ blocks corresponding to
\eqn\decomp{H^*(X)=\sigma H^*(B)\oplus \pi^* H^*(B)}

A major part of the discussion will be concerned with the question
whether one has to transform just
the Chern characters $ch(V)$ and $ch(\tilde{V})$ themselves
or to use appropriately twisted charges $Q$, $\tilde{Q}$
(like in the usual Fourier-Mukai transformation where the charges
have twists with
$\sqrt{Td(X)}$) and similarly the related question whether the
transformation process itself has to be somewhat twisted.

\subsec{Description of the bundle by spectral cover data}
Let us first recall the spectral cover description.
The $SU(n)$ bundle $V$ on $X$ of $c_1(V)=0$ decomposes on the typical
fibre $E$ (where $V$ is assumed to be semistable) as a sum
$\oplus_i {\cal L}_i$
of line
bundles of degree zero and each of the ${\cal L}_i$
corresponds\foot{The double interpretation of $E$ as pointset $E_1$
resp. parameter space $E_2$ for degree zero line bundles on $E_1$ is
formalized by introducing the Poincar\'e bundle ${\cal P}$ on $e_1\x E_2$
which restricts on $E_1\x Q$ to  ${\cal L}_Q$; actually one uses the
symmetrized version ${\cal P}={\cal O}(\Delta - p\x E_2 - E_1\x p)$.}
(having chosen the distinguished reference point $p$ as origin) to a
welldefined point $Q_i$ on $E$ (these points sum up to zero as
$det(V)=1$). When the reference point is globalized by the section
$\sigma$ the variation of the $n$ points in a fibre lead to a
hypersurface $i: C\hookrightarrow X$, a ramified $n$-fold cover (the 'spectral
cover') of $B$. The equation $s=0$ of $C$ involving the section $s$ of
${\cal O}(\sigma)^n$ can in the process of globalization still be
twisted by a line bundle ${\cal M}$ over $B$ of $c_1({\cal M})=\eta$,
i.e. $S$ can be actually a section of
${\cal O}(\sigma)^n\otimes {\cal M}$
and the cohomology class of $C$ in $X$ is given by
$$C= n\sigma + \eta$$
Now $V$ will be induced as $V=p_* {\cal R}$
from a line bundle ${\cal R}$ over the
$n$-fold cover $p:X\x _B C\rightarrow X$, i.e. generically will the
fibre of $V$ over a point $x\in X$ with the $N$ preimages
$\tilde{x}_i$ be given by the sum of the fibre of ${\cal R}$
at the $\tilde{x}_i$. If one takes for ${\cal R}$ the
global version of
${\cal P}$ on has indeed that {\it fibrewise}
$V=p_* {\cal P}$ as
$p_*$ sums up the line bundles which makes ${\cal P}$ out of the
points collected in the fibre $C_b$ of $C$ over $b\in B$, which
themselves corresponded to the line bundles summands of $V$ on
$E_b$. As the twist by a line bundle $L$ over $C$ leaves the fibrewise
isomorphism class unchanged ($L$ being locally trivial along $C$) the
construction generalizes to
\eqn\spectrcov{V=p_*(p_C^*L\otimes {\cal P})}
where $p$ and $p_C$ are the projections on the first and second factor
of $X\times_B C$
$$\matrix{\;\;\; X\times_B C&\buildrel p_C \over \longrightarrow&C\cr
\scriptstyle{p} \biggr \downarrow & &
\scriptstyle{\pi_C}\biggr \downarrow \cr
     \;\;\; X&\buildrel {\pi_1} \over\longrightarrow&B}$$

The condition $c_1(V)=0$ translates to a
fixing\foot{$\pi_*(c_1(L))=-\pi_*(c_1(C)-c_1)/2)$}
of $c_1(L)$ in $H^{1,1}(C)$ up to a class in $ker \;
\pi_*:H^{1,1}(C)\rightarrow H^{1,1}(B)$; such a class is known to be
of the form $\gamma=\lambda (n\sigma - (\eta - nc_1))$ with $\lambda$
half-integral.

\subsec{FM transform}
For the description of the FM transform we will
instead of working on $X\times_B C$ work on $X\times_B \tilde X$
\eqn\mats{\matrix{\;\;\; X\times_B{\tilde X}&\buildrel p_2 \over
\longrightarrow&{\;\;\tilde X}\cr
\scriptstyle{p_1} \biggr \downarrow & &
\scriptstyle{\pi_2}\biggr \downarrow \cr
     \;\;\; X&\buildrel {\pi_1} \over\longrightarrow&\;\;B}}
where $\tilde X$ is the compactified relative Jacobian of $X$. $\tilde
X$ parameterizes torsion-free rank 1 and degree zero sheaves of the
fibres of $X\to B$  and it is actually isomorphic with $X$  (see
\bart\ or \RPo). We will then identify $\tilde X$ and $X$.

The bundle $V$ is then given by \foot{In what sequel we will identify a
bundle with the locally free sheaf of its sections, so that the terms
bundle and locally free sheaf are used interchangeable}
\eqn\bundv{V=R^0p_{1*}(p_2^*(i_*L)\otimes {\cal P})}
where
\eqn\poincdef{{\cal P}=
{\cal O}(\Delta - \sigma \x \tilde{X} - X\x \tilde{\sigma}- c_1(B))}
is the Poincar\'e sheaf normalized to
make ${\cal P}$ trivial along $\sigma \x \tilde{X}$
and $X\x \tilde{\sigma}$\foot{Here we denote by ${\cal O}(\Delta)$ the
dual of the ideal sheaf of the diagonal. Neither ${\cal O}(\Delta)$
nor $\cal P$ are line bundles due to the presence of singular fibres,
but they are torsion-free and rank 1}.

Now let us determine the FM-transform.
For this we make use of the fact that the representation of $V$ by the
$(C,L)$ data already looks in itself like a FM transform; so when we
want to describe now the FM transform of $V$ this is practically the
double transform of $i_* L$; but it is known that the inverse
transform of FM is not precisely FM itself again but a slightly
twisted version of the first transform as we shall see later in this
section; so only this twisted version would bring us back
from
$V$ to $i_* L$, or said differently, if we now start on $V$ by making
the original FM transform we will get
$i_* L$ times the inverse twist.

So with
$V=p_{1*}(p_2^*(i_*L)\otimes {\cal P})$
the 'fibrewise dual' bundle is given\foot{As mentioned above
the logic of the procedure is that
this follows from the fact that the inverse FM transform is given by
$W=R^1p_{2*}(p_1^*V\otimes {\cal P}^*\otimes
p_2^*\pi_2^*K_B^{-1} )
=R^1p_{2*}(p_1^*V\otimes {\cal P}^*)\otimes \pi_2^*K_B^{-1}$
so that we have (with $W=i_* L$)
$W=\tilde V \otimes \pi_2^*K_B^{-1}$.}
in terms
of the spectral cover data by\foot{
This is fibrewise just the fact that
$\tilde{V}=\oplus_i (\pi_2)_*(\pi_1^* ({\cal L}_{Q_i})\otimes {\cal
P}^*)=\oplus_i {\cal O}_{Q_i}=i_*1_L$
with the trivial line bundle $1_L$ on $C=\cup \{ Q_i \}$,
for note that
${(\pi_1^* ({\cal L}_{Q_i})\otimes {\cal P}^*)|_{E_1\x q}}$
$=\atop \neq$ ${{\cal O}_E}$ for $q$$=\atop \neq $$Q_i$.}

\eqn\dulav{\eqalign{{\tilde V}&=R^1p_{2*}(p_1^*(V)\otimes
 {\cal P}^*)\cr
                              &=i_*L\otimes \pi_2^*K_B}}

Note that one wants now to show that there exists a matrix $M$ in the block
diagonal form as above in \precisematrix\
which relates $V$ to $\tilde{V}$ or $V$ to $i_* L$
(this is not a big difference because the latter option
is simply the inverse process as $V$ is the dual of $i_* L$).

Before going on we would like to recall some well-known facts
about the FM transform for elliptic fibrations (cf. \bart\ ,
\bridge\ or \RPo) that explain the facts mentioned above.

We define the Fourier-Mukai functors $S^i$, $i=0,1$
by associating with every sheaf $V$ on $X$ the sheaf
$S^i(F)$ on
$X$ (where $X$ and $\tilde X$ are identified)
\eqn\sif{S^i(V)=R^ip_{1*}(p_2^*(V)\otimes {\cal P})}
where ${\cal P}$ denotes the Poincar\'e sheaf \poincdef\ on the fibre product.
It can be also described as (cf. \FMW\ )
\eqn\poinc{{\cal P}={\cal I}^*\otimes p_1^*{\cal O}(-\sigma)\otimes
p_2^*{\cal O}(-\sigma)\otimes q^*K_B}
with $q=\pi.p_1=\pi.p_2$ and ${\cal I}={\cal O}(\Delta)^*$
the ideal sheaf of the diagonal
immersion $\delta: X \rightarrow X\times_B X$.

We can also define the inverse Fourier-Mukai functors $\hat{S}^i$, $i=0,1$
by associating with every sheaf $V$ on X the sheaf
\eqn\infm{\hat{S}^i(V)=R^ip_{2*}(p_1^*(V)\otimes {\cal
P}^*\otimes
q^*K_B^{-1})}
The relationship between these functors is more neatly stated if we consider
the associated functors between the derived categories of complexes of
coherent sheaves bounded from above.
\eqn\assfd{\eqalign{S&: D^{-}(X)\rightarrow D^{-}(X);\ \ S({\cal
G})=Rp_{1*}(p^*_2({\cal G})
\otimes {\cal P})\cr
\hat{S}&: D^{-}(X)\rightarrow D^{-}(X);\ \ \hat{S}({\cal
G})=Rp_{2*}(p^*_1({\cal G})
\otimes {\cal P}^*\otimes q^*K_B^{-1})}}
Proceeding as in Theorem 3.2 of \bart\ and taking into account (\RPo, Lemma
2.6) one obtains\foot{For any complex
${\cal G}$ (or any element in the category) with cohomology
sheaves ${\cal G}^i$, the cohomology sheaves of the shifted complex
${\cal G}[n]$ are  ${\cal G}[n]^i={\cal G}^{i+n}$}
an invertibility result\foot{Concerning the meaning of the -1
shift, consider a complex given by a single sheaf
$V$ located at the ``degree zero'' position.
$\hat{S}({S}(V))= V[-1]$ means that the complex $\hat{S}({S}(V))$
has only one cohomology sheaf, which is $V$, but located at ``degree
1'', $[\hat{S}({S}(V))]^1=V,\ \ [\hat{S}({S}(V))]^i=0, i\neq 1$.
When $S^0(V)=0$ the
complex $S(V)$ reduces to a single sheaf, which is
the unique FM transform $S^1(V)$, but located at ``degree 1'', that
is,
$S(V)=S^1(V)[-1]$ and the
complex $\hat{S}(S(V))=\hat{S}(S^1(V))[-1]$ has two cohomology
sheaves, one at degree 1 which is $\hat{S}^0(S^1(V))$, and
one at degree 2 which is $\hat{S}^1(S^1(V))$. So one has
$\hat{S}^0(S^1(V))=V, \ \ \hat{S}^1(S^1(V))=0$.}:
\eqn\fcd{S(\hat{S}({\cal G}))= {\cal G}[-1], \ \ \hat{S}({S}({\cal G}))= {\cal
G}[-1]}
{\it Remark:}
Even when ${\cal G}$ is a single sheaf $V$, that is, a complex
with
$V$ at degree 0 and no other terms, $S(V)$ is an object of the derived
category, or
a complex whose cohomology sheaves are the sheaf FM transforms $S^i(V)$.
It is then interesting to know when only one of the FM sheaf
transforms is different from zero, the so-called WIT$_i$ condition. This
condition is better studied when the sheaf $V$ is flat over $B$
(for instance, a vector
bundle or a sheaf $i_*L$ where $L$ is a line bundle on a spectral cover
$C$ flat of degree $n$ over $B$\foot{We mean that $\pi|_C:C\ra B$ is a flat
morphism of degree
$n$, that is, all its fibres consist of $n$ points (mayby counted more than
once); $i:C\ra X$ is the embedding.}). In this case, the vanishing of
$S^1(V)$ or WIT$_0$ condition, is equivalent to WIT$_0$ on every fibre,
that is $S^1(V)=0\iff  H^1(X_t,V_t\otimes
{\cal L}_x)=0$ for every point $x$
where $X_t$ is the fibre over
$t=\pi_2(x)$, $F_t=F|_{X_t}$ and ${\cal L}_x$
is the rank one torsion free sheaf of degree 0 on $X$ defined by $x$;
a sheaf supported fibrewise by points, as $i_*L$, is always
WIT$_0$.
The WIT$_1$ condition ($S^0(V)=0$) is not a fibrewise condition. If
$V$ is WIT$_1$ on every fibre, (that is, $H^0(X_t,V_t\otimes
{\cal L}_x)=0$ for every point $x$) then $V$ is globally
WIT$_1$ and the only FM transform $S^1(V)$ is flat over $B$ (\RPo
2.11). This happens for a vector bundle
$V$ fibrewise semistable of degree zero.
Conversely, if $V$ is
WIT$_1$ the flatness of $S^1(V)$ is necessary to ensure that
$V$ is fibrewise WIT$_1$. A typical example is a rank $n$ vector bundle
$V$ fibrewise of degree zero, which is only semistable on the generic
fibre. It is still WIT$_1$ but fails to be so at those fibres where it
is unstable; this reflects the fact that the spectral cover $C$
contains those fibres so that neither $C$ nor $S^1(V)=i_*L$ are flat
over $B$.

Later we will assume $c_1(V)=0$ or at least that $V$ has degree 0 and is
semistable when restricted to the  fibre; this
is the most important case to our purposes and the bundles given by the
spectral
cover construction are of this kind.
\subsec{K3 case}
Let us consider first the two-dimensional case of $X=K3$ and assume here
the decomposition
\eqn\decompok{\eqalign{
H^0(X)&={\bf C}\cr
\oplus&\cr
H^2(X)&={\bf C}\sigma \oplus {\bf C}\cr
\oplus&\cr
H^4(X)&={\bf C}}}
The class
of the spectral curve $C$ on which $i_*L$ is supported is $C=n\sigma+kF$.

The Chern characters of $i_*L$ can be obtained
using Grothendieck-Riemann-Roch for the embedding $i:C\rightarrow
{\tilde X}$
\eqn\grrfe{ch(i_*L)Td({\tilde X})=i_*(ch(L)Td(C))}
one gets
\eqn\grrd{ch(i_*L)=(0,C,n)}
further $ch(V)$ ($V$ is the only FM transform of $i_*L$) is given by
\eqn\chv{ch(V)=(n,0,-k)}
Now let us introduce
the new functor
$T(\cdot)=S (\cdot )\otimes \pi^* K_B^{-1/2}
=S (\cdot )\otimes {\cal O}(F)$ (and similarly $T^i(\cdot)=S^i (\cdot )\otimes
\pi^* K_B^{-1/2} =S^i(\cdot )\otimes {\cal O}(F)$) so that $T^0(i_*
L)=V\otimes {\cal O}(F)$.  We get\foot{This holds
more generally: if only
$T^i({\cal F})$ is non-zero, one has
$(-1)^ich(T^i({\cal F}))=M\cdot ch(i_*L)$
and the same formula for $\hat{T}$ (\bart, eqn.(4.1)). As
the only inverse transform of $T^i({\cal F})$ is $\hat{T}^{1-i}$, one has
$ch(\hat{T}^{1-i}(T^i({\cal F})))=-M^2\cdot ch({\cal F})$,
consistent with $\hat{T}^{1-i}(T^i({\cal F}))=F$ and $M^2=- {\bf
id}_4$.}
\eqn\tfuc{\eqalign{ch(T^0(i_* L))&=ch(V)(1+F)\cr
                        &=M\cdot ch(i_*L)}}
where we reach the matrix we wanted
$$
M=\pmatrix{0&1&0&0\cr
-1&0&0&0\cr
0&0&0&1\cr
0&0&-1&0}$$
The functor $T$ was introduced
\bart\ because its inverse transform in the sense of \fcd\ is the ``natural''
one: $T$ is the FM transform w.r.t. the sheaf $\bar
{\cal P}={\cal P}\otimes q^*{\cal O}(F)$  and its inverse functor $\hat{T}$
is the FM transform w.r.t. the dual sheaf $\bar {\cal P}^*$.
$S$ does not have this property as its inverse transform is not the
FM transform with respect to ${\cal P}^*$, but this
twisted by $q^*K_B^{-1}$. So we just
divided this up in two parts of
$K_B^{-1/2}$, distributed among the original and the
inverse transform.

\subsec{Calabi-Yau threefold case}
So far we considered FM transformation of $V$ with vanishing first
Chern class in the context of \FMW\ . Let us now
describe the topological invariants of the relative FM transform for
a coherent sheaf ${\cal F}$ on an elliptic Calabi-Yau threefold. Therefore we
start again from \mats.
Applying GRR for $p_1$ we get (${\cal G}$
can be an object of the derived category)
\eqn\grrfp{ch(S({\cal G}))=p_{1*}[p^*_2({\cal G})\cdot ch({\cal P}) \cdot
Td(T_{X/B})]}
where $Td(T_{X/B})=1-{1\over 2}c_1+{1 \over 12}(13c_1^2+12\sigma
c_1)- {1\over 2}\sigma c_1^2$ (with  $c_1=\pi^*c_1(B)$)
denotes the Todd class of the relative tangent bundle $T_{X/B}=T_X/\pi^*T_B$.
Note that  $S({\cal G})$ is a complex (or an object of
the derived category) and then its Chern character is
\eqn\der{ch(S({\cal G}))=\sum_i(-1)^ich(S^i({\cal G}))}
To compute \grrfp\ note first that $ch({\cal I})=1-ch(\delta_*{\cal O}_X)$
with the diagonal immersion $\delta$. Riemann-Roch gives
\eqn\srr{ch(\delta_*{\cal O}_X)Td(X\times_B X)
=\delta_*(ch({\cal O}_X)Td(X))}
where one has the expressions for $Td(X)$ and $Td(X\times_B X)$ given by
\eqn\tdx{\eqalign{Td(X)&= 1+{1\over 12}(c_2+11c_1^2+12\sigma c_1)\cr
        Td(X\times_B X)&=p_2^*Td(X)p_1^*Td(T_{X/B})}}
The Chern character of the ideal sheaf is then given by
(with the diagonal class $\Delta=\delta_*(1)$)
\eqn\chois{\eqalign{ch({\cal I})&=1-\delta_*(1)-{1\over
2}\delta_*(c_1)+\delta_*(
\sigma\cdot c_1)+{5\over 6}\delta_*(c_1^2)+{1\over 2}
\delta_*(\sigma c_1^2)\cr
&=1-\Delta-{1\over 2}\Delta\cdot p_2^*c_1+\Delta\cdot p_2^*(
\sigma\cdot c_1)+{5\over 6}\Delta\cdot p_2^*(c_1^2)+{1\over
2}
\Delta\cdot p_2^*(\sigma c_1^2)}}
Defining the numerical invariants (with $F$ the elliptic fibre class)
\eqn\numin{\eqalign{n&=rk{{\cal G}}, \ \ \ \ s=ch_3({\cal G})\cr
                    d&=ch_1({\cal G})\cdot F, \ \ \ \
                    g=ch_1({\cal G})\cdot\sigma\cdot c_1\cr
                    c&=ch_2({\cal G})\cdot\sigma \ \ \ \
                    f=ch_2({\cal G})\cdot c_1}}
we get for the Chern characters of $S({\cal G})$
\eqn\chS{\eqalign{ch_0(S({\cal G}))&=d\cr
  ch_1(S({\cal G}))&=ch_1({\cal G})-(d+n)\sigma-p_{1*}(p_2^*ch_1({\cal
G})\sigma)+
  p_{1*}(p_2^*ch_2({\cal G}))-{3\over 2}dc_1\cr
  ch_2(S({\cal G}))&=ch_2({\cal G})-2ch_1({\cal G})(c_1+\sigma)+\sigma p_{1*}
 (p_2^*(ch_1({\cal G})\sigma))-\sigma p_{1*}(p_2^*ch_2({\cal G}))\cr
&\ \ \  +{25\over 12}dc_1^2+(s+2g-c-
{3\over 2}f)F+({1\over 2}n+d)\sigma c_1\cr
ch_3(S({\cal G}))&=-({1\over 6}n\sigma c_1^2-{1\over 2}\sigma
c_1^2d-{1\over 2}g+f+c)}}
Note that if there is only one non-vanishing transform $S^i$, its Chern
character is computed from \chS\ by
$ch(S^i({\cal G}))=(-1)^ich(S({\cal G}))$ due to \der.

Similar calculations can be done for the inverse FM transform.
\eqn\chSS{\eqalign{ch_0(\hat{S}({\cal G}))&=d\cr
  ch_1(\hat{S}({\cal G}))&=-ch_1({\cal
G})+(d-n)\sigma+p_{2*}(p_1^*ch_1({\cal G})\sigma)+
  p_{2*}(p_1^*ch_2({\cal G}))+{3\over 2}dc_1\cr
  ch_2(\hat{S}({\cal G}))&=-ch_2({\cal G})
-2ch_1({\cal G})(c_1+\sigma)+\sigma p_{2*}
 (p_1^*(ch_1({\cal G})\sigma))+\sigma p_{2*}(p_1^*ch_2({\cal G}))\cr
&\ \ \  +{25\over 12}d c_1^2+(s+2g+c+
{3\over 2}f)F+(d-{1\over 2}n)\sigma c_1\cr
ch_3(\hat{S}({\cal G}))&=-({1\over 6}n\sigma c_1^2-{1\over
2}\sigma c_1^2d+ {1\over 2}g+f+c)}}

If we consider now a sheaf $V$  and write its Chern
character as
\eqn\chV{ch_0(V)=n,\
ch_1(V)=x\sigma+S,\
ch_2(V)=\sigma \eta+aF,\
ch_3(V)=s}
($\eta, S\in p_2^*H^2(B)$), then by \chS\ and \chSS, the Chern character of the
FM of $V$ and of the inverse FM of
$V$ are
\eqn\chVS{\eqalign{ch_0(S(V))&=x\cr
ch_1(S(V))&=-n\sigma+\eta-{1\over 2}xc_1\cr
ch_2(S(V))&=({1\over 2}nc_1-S)\sigma+(s-{1\over 2}\eta c_1\sigma
+{1\over12}xc_1^2\sigma)F\cr
ch_3(S(V))&=-{1\over 6}n\sigma c_1^2-a+{1\over 2}\sigma c_1
S}}
and
\eqn\chVSS{\eqalign{ch_0(\hat{S}(V))&=x\cr
ch_1(\hat{S}(V))&=-n\sigma+\eta+{1\over 2}xc_1\cr
ch_2(\hat{S}(V))&=(-{1\over 2}nc_1-S)\sigma+(s+{1\over 2}\eta c_1\sigma
+{1\over12}xc_1^2\sigma)F\cr
ch_3(\hat{S}(V))&=-{1\over 6}n\sigma c_1^2-a-{1\over 2}\sigma c_1
S+x\sigma c_1^2}}
Using the decompostion of the cohomology we find
$$
ch(V)=\pmatrix{n\cr x\cr S\cr \eta \cr a \cr s}\ ,\qquad
ch(\hat{S}(V))=\pmatrix{0\cr -n\cr \eta+{1\over 2}xc_1\cr -{1\over 2} n
c_1-S\cr
s+{1\over 2}\eta c_1\sigma+{1\over12}xc_1^2\sigma\cr
-{1\over 6}n\sigma c_1^2-a-{1\over 2}\sigma c_1
S+x\sigma c_1^2}
$$
If we multiply $ch(\hat{S}(V)$ by the Todd class
$Td(N)=1-{1\over2}c_1+{1\over12}c_1^2$,\foot{We always confuse the
normal bundle $N$ with its pull-back $\pi^* N$ to $X$} we get
$$
ch(\hat{S}(V)\cdot Td(N)=\pmatrix{x \cr
-n \cr \eta \cr -S \cr s-{1\over 12}x\sigma c_1^2 \cr
-a+x\sigma c_1^2
}
$$
When $V$ is fibrewise of degree 0, which is the case we are mainly
interested in (bundles constructed from spectral covers are of this
kind), we have $x=0$ and then
\eqn\Mmat{
\pmatrix{x \cr
-n \cr \eta \cr -S \cr s \cr
-a
}=\pmatrix{0&1&0&0&0&0\cr
-1&0&0&0&0&0\cr
0&0&0&1&0&0\cr
0&0&-1&0&0&0\cr
0&0&0&0&0&1\cr
0&0&0&0&-1&0}\cdot\pmatrix{n\cr x\cr S\cr \eta \cr a \cr s}
}
that is
\eqn\MmathS{Td(N)\cdot ch(\hat{S}(V))=M\cdot ch(V)\,.}
So for fibrewise
degree 0 and semistable bundles $V$ we have  $\hat{S}^0(V)=0$
and
$\hat{S}^1(V)=i_*L$, so that \Mmat\ is equivalent to
$$
Td(N)\cdot ch(i_*L)=Td(N)\cdot ch(\hat{S}^1(V))=-M\cdot ch(V)\,.
$$

There is an analogous equation of \MmathS\ for the direct FM
transform $S(V)$; proceeding as above one shows that if $x=0$
then
\eqn\MmatS{Td(N^{-1})\cdot ch(S(V))=M\cdot ch(V)\,.}
For
fibrewise degree 0 and semistable bundles $V$ we have
$S^0(V)=0$ and if we write
$S^1(V)=i_*\bar{L}$, we have
$$
Td(N^{-1})\cdot ch(i_*\bar{L})=Td(N^{-1})\cdot ch(S^1(V))=-M\cdot
ch(V)\,.
$$

How have these results to be interpreted ?
We will show in the appendix why
in the three-fold case
the sole use
of the T-functor known from the $K3$-case
to give the map between the Chern classes of the bundle and its dual
is insufficient to exhibit as transformation matrix
the adiabatic extension \precisematrix\
of the usual T-duality matrix \ellmatrix\
on the fibre.

Rather one has to invoke the precise definition of the twisted charge
relevant here. Recall from the usual (full, not fibrewise)
Fourier-Mukai transform that one has actually to take the Mukai vector
which has a twist by $\sqrt{Td(X)}$ and not just the Chern character.

Now there are two possible alternative routes of procedure: one can
either twist the operation a little bit and stick to the usual
Mukai-vector (this is described in the next subsection) or keep the
operation and make the twist in the usual Mukai-vector more 'relative'
as described in what follows.

Similarly to the usual twist in the Mukai-vector
here in the fibrewise situation a twist by $\sqrt{Td(N)}$
with the normal bundle $N=j^*TX/TB=j^*T_{X/B}$ plays a role
(where $j:B\hookrightarrow X$).
Now suppose you have a complex of sheaves ${\cal G}$ on $X$ and
you twist the standard definition of $ch$ to define  a ``charge''
\eqn\ch{Q({\cal G})=\sum(-1)^ich({\cal G}^i)(\sqrt{Td(N)})^{|i+1|}}
where ${\cal G}^i$ are the cohomology sheaves of the complex ${\cal G}$ and
$|i+1|=(-1)^{i+1}$.

Then
\eqn\chd{\eqalign{Q(\hat{S}(V))&=-ch(\hat{S}^1(V))(\sqrt{Td(N)})\cr
                  Q(V)&=ch(V)(\sqrt{Td(N)})^{-1}}}
and thus\foot{For a single sheaf $V$, we are writting $Q(V)$ in the sense
that $V$ is understood as a complex with $V$ at degree 0 and no other terms.}
\eqn\mat{Q(\hat{S}(V))=M\cdot Q(V)}
Since $\hat{S}^0(V)=0$ and $\hat{S}^1(V)=i_*L$ we have
$\hat{S}(V)=i_*L[-1]$ as complexes and then \chd\ is equivalent to
$$
Q(i_*L[-1])=M\cdot Q(V)=M\cdot Q(S(i_*L))
$$

This represents one way to arrange the quantities
involved to get the M matrix. In view of the internal twist by $Td(T_{X/B})$
in the FM transform \grrfp\  the charge definition \ch\  may not come
as a surprise\foot{note the close relation $j^*T_{X/B}=N$ between
$T_{X/B}=TX/\pi^* TB$ and $N=(j^* TX)/TB$ where
$j:B\hookrightarrow X$; cf. also the motivation for the T-functor
in the $K3$ case
described at the end of that subsection and the appendix for the role of
$Td(N)$}. As outlined above still a different route can be taken. This
is described in what follows. 

\subsec{f-map}

The effect of the FM transform in cohomology is usually
described by means of the so-called ${\bf f}$-map. For a
fibrewise FM transform, we can also introduce a relative
version
${\bf f}_r: H(X)\to H(X)$ of the
${\bf f}$-map. It is defined by
$$
{\bf f}_r(x)=p_{1*}(p_2^*(x)\cdot Z_r)
$$
where $Z_r=\sqrt{p_2^*Td(T_{X/B})}\cdot ch({\cal P})\cdot
\sqrt{p_1^*Td(T_{X/B})}$. Then
$$
\sqrt{Td(T_{X/B})}\cdot ch(S(V))={\bf f}_r\big(ch(V)\cdot
\sqrt{Td(T_{X/B})}\ \big)
$$
If we consider instead of $ch(V)$ the effective charge
given by the Mukai vector
\eqn\mukf{Q(V)=ch(V)\cdot\sqrt{Td(X)}}
then the effect of ${\bf f}_r$ on $Q(V)$ is described by
\eqn\fmapr{{\bf f}_r(Q(V))=Q(S(V))}
But if we modify the definition of ${\bf f}_r$ to ${\bf f}: H(X)\to
H(X)$,
\eqn\fmapdef{{\bf f}(x)=p_{1*}(p_2^*(x)\cdot Z)}
with $Z=\sqrt{p_2^*Td(X)}\cdot ch({\cal P})\cdot
\sqrt{p_1^*Td(X)}$, then by \MmatS, the effective charge of $V$
transforms to
\eqn\fmap{{\bf f}(Q(S(V)))=M\cdot Q(V)}
when $x=0$.  Then if $V$ is moreover fibrewise semistable so that
$S^0(V)=0$ and
$S^1(V)=i_*L$, then
$$
-{\bf f}(Q(i_*L))=-{\bf f}(Q(S^1(V)))={\bf f}(Q(S(V)))=M\cdot Q(V)
$$

\subsec{Fibrewise T-duality on D-branes at the sheaf level}

Our next goal is to describe T-duality on the $T^2$ fiber
maps given in \Dtual\ at the sheaf level.

Let us consider the skyscraper
sheaf ${\bf C}(x)$ at a point $x$ of $X$. It is a WIT$_0$ sheaf and its FM
transform $S^0({\bf C}(x))$ is a torsion-free rank one sheaf $L_x$ on the
fibre of $X$ over $\pi(x)$ \foot{With the identification $X\simeq\tilde X$
the point
$x$ corresponds precisely to $L_x$ (see \bart\ or \RPo\ )} as we expect
from \Dtual\
and thus we see $D0\rightarrow D2_F$.

For the topological invariants we have indeed $n=x=a=0, S=\eta=0, s=1$ and then
\eqn\fmap{ch_i(S^0({\bf C}(x)))=0, \ \ i=0,1,3, \ \ \ \ ch_2(S^0({\bf
C}(x)))=F}

If we start with ${\cal F}={\cal O}_\sigma$; proceeding as in (3.16) of
\RPo\ we
have\foot{Our formulas differ from those in \RPo\ because we are using a
different
Poincar\'e sheaf}
\eqn\fmsigma{\eqalign{S^0({\cal O}_\sigma)&={\cal O}_X\,,\qquad
S^1({\cal O}_\sigma)=0\cr S^0({\cal O}_X)&=0\,,\quad\qquad
S^1({\cal O}_X)={\cal O}_\sigma\otimes \pi^*K_B}}
Then ${\cal O}_\sigma$ transforms to the structure sheaf of $X$ and ${\cal
O}_X$ transforms to a line bundle on $\sigma$ as we expect from
\Dtual
since $D4_B\leftrightarrow D6$. We have as before the transformations at the
cohomology level;
\eqn\bas{n=0, \ \ x=1, \ \ S=0, \ \ \eta={
1\over 2}c_1, \ \
a= 0, \ \ s={1\over 6}\sigma c_1^2}
and then we get
\eqn\resu{ch_0(S^0({\cal O}_\sigma))=1,\quad ch_i(S^0({\cal
O}_\sigma))=0,\ \ i=1,2,3}
Finally, let us consider a sheaf ${\cal F}$ on $B$; by \fmsigma\ we have
\eqn\fmsigmatwo{\eqalign{&S^0({\cal O}_\sigma\otimes \pi^*{\cal
F})=\pi^*{\cal F}\,,\qquad S^1({\cal O}_\sigma\otimes
\pi^*{\cal F})=0\cr &S^0(\pi^*{\cal F})=0\,,\phantom{{\cal O}_\sigma\otimes
\pi^*{\cal F}}\ \quad S^1(\pi^*{\cal F})={\cal O}_\sigma\otimes \pi^*{\cal
F}\otimes
\pi^*K_B}}
Then, a sheaf ${\cal
O}_\sigma\otimes \pi^*{\cal F}=j_*{\cal F}$ ($j: B\ra X$ is the section)
supported on a curve
$\tilde C$ in $B$ embedded in $X$ via $j$ transforms to a  sheaf on the
elliptic surface supported on the inverse image of $\tilde C$ in $X$ and
vice versa.
This is what we expected form the map $D2_B\leftrightarrow D4_F$ of \Dtual.

Then, even at the sheaf level we have the
relations
\Dtual\ appropriate for the fibrewise T-duality on D-branes
\eqn\Dtual{\eqalign{D4_B&\rightarrow \tilde{D6}\cr
D2_B&\rightarrow \tilde{D4}_{\tilde{F}}\cr
D0&\rightarrow \tilde{D2}_{\tilde{F}}}}

\newsec{Period monodromies and bundle automorphisms}
Kontsevich proposed to consider mirror symmetry as an categorical
equivalence between the bounded derived category $D(X)$ of $X$ and
Fukaya's $A_{\infty}$ category of Lagrangian submanifolds for the mirror
$\hat{X}$. The object of the derived category $D(X)$ is a complex of
coherent sheafs on $X$ resp. the object of Fukaya's $A_{\infty}$ category
is a Lagrangian submanifold with flat $U(1)$ bundles on it.

It is now expected that under the proposed categorical equivalence between
$D(X)$ and $A_{\infty}$ the monodromy of the SLAG 3-cycle is mapped to
certain automorphisms in $D(X)$ (Fourier-Mukai transformations). We
are going to make these relations explicit and check them thereby for
two models\foot{which are just $K3$ fibrations, not elliptic ones, so
somewhat outside our main applications}.
\subsec{The models}
So we consider the models
${\bf P}_{1,1,2,2,2}^4[8]$ and ${\bf P}_{1,1,2,2,6}^4[12]$
given by degree 8 resp. 12 hypersurfaces in the respective weighted
projective spaces. Among these hypersurfaces are
\eqn\clasgeo{z_1^8+z_2^8+z_3^4+z_4^4+z_5^4=0}
and
\eqn\clasgeotw{z_1^{12}+z_2^{12}+z_3^6+z_4^6+z_5^2=0}
These models were studied in \dave, \hoskl. At $z_1=z_2=0$ both models have
a curve $C$ (of genus three resp. two) of $A_1$ singularities,
which when blown up leads to an exceptional
divisor $E$ in $X$.
$E$ has the structure of a ruled surface ($P^1$ fibration) over $C$;
let $l$ denote its fibre.
A second divisor class $L$ (besides $E$) comes from the fact that both
models are $K3$ fibrations over $P^1$
(the degree one polynomials
generate a linear system
$|L|$ which projects $X$ to ${\bf P}^1$ with fiber $L=K3$).
$L$ and $E$ together generate $H_4(X,{\bf Z})$.
The generators of the complexified K\"ahler cone in $H^2(X,{\bf Z})$
are chosen to be $(E,L)$, thus the generic
K\"ahler class can be written as $t=t_1E+t_2L$ where $(t_1,t_2)$
are coordinates on the K\"ahler moduli space of $X$.

A second linear system $|H|$ is generated by degree two polynomials
and related to the first by $|H|=|2L+E|$ (for both models).
Let $4h$ denote the intersection of two general members
of $H$ and $L$ in the degree 8 model, which is a plane quartic
(think of a hyperplane section of the $K3=P^3[4]$)
which can be specialized to a sum of four lines.
Similarly one can define $2h=H\cdot L$ in the degree 12 model.
$h={1\over 4}H\cdot L$ and $l={1\over 4}H\cdot E$
in the degree 8 model resp.
$h={1\over 2}H\cdot L$ and $l={1\over 2}H\cdot E$
in the degree 12 model
generate $H_2(X,{\bf Z})$ in both cases. They are dual to $H$ and $L$ in
the sense that the relations
$L\cdot l=1, \ \ L\cdot h=0, \ \ H\cdot l=0, \ \ H\cdot h=1$
hold (which mean for $E$ that $E\cdot l=-2,\ \ E\cdot h=1$).

Let us recall now the intersection relations for both models.
As they are $K3$ fibrations one has $L^2=0$. The remaining
intersections are as follows.
For ${\bf P}^4_{1,1,2,2,2}[8]$ one has\foot{one has also the relations
$H^3=8, \ \ H^2\cdot L=4$ and $H\cdot E\cdot L=4$}
\eqn\scoh{\eqalign{& E^3=-16, \ \ E^2\cdot L=4\cr
                   & c_2(X)\cdot L=24, \ \ c_2(X)\cdot E=8}}
and for ${\bf P}^4_{1,1,2,2,6}[12]$ one has\foot{one has also the
relations
$H^3=4, \ \ H^2\cdot L=2$ and $H\cdot E\cdot L=2$}
\eqn\scoh{\eqalign{& E^3=-8, \ \ E^2\cdot L=2\cr
                   & c_2(X)\cdot L=24, \ \ c_2(X)\cdot E=4}}
%

The mirror family $\hat X$ of $X$ can be obtained by applying
the Greene-Plesser construction \grepl,
which is given by ${p=0}/G$ with
\eqn\mir{p=z_1^8+z_2^8+z_3^4+z_4^4+z_5^4
-8\psi z_1z_2z_3z_4z_5-2\phi z_1^4z_2^4=0}
for ${\bf P}^4_{1,1,2,2,2}[8]$ with $G={\bf Z}_4^3$
and for ${\bf P}^4_{1,1,2,2,6}[12]$
\eqn\mirtw{p=z_1^{12}+z_2^{12}+z_3^6+z_4^6+z_5^2
-12\psi z_1z_2z_3z_4z_5-2\phi z_1^6z_2^6=0}
with $G={\bf Z}_6^2\times {\bf Z}_2$. Note that $\phi$ and $\psi$
parameterize the moduli
space of complex structures of the mirror $\hat X$
which does not get any $\alpha'$ corrections.

Of course the complex structure moduli space of the mirror is
identified with the K\"ahler moduli space of the original Calabi-Yau
and under this identification monodromies on the K\"ahler side are replaced
by there transpose-inverse. 
The $N=2$ prepotential ${\cal F}$ determines
the vector of periods of the holomorphic three-form
$\hat\Omega$ on the mirror manifold $\hat X$ which is given
in the $(E,L)$ basis by (note the identification just mentioned)
$$\Pi=\pmatrix{2{\cal F}-t^iF_i\cr
                F^1-2F^2\cr
                F^2\cr
                1\cr
                t_1\cr
                t_2}$$
Note that the order\foot{
The 2. entry is not $F^1$ as one
comes from the $H,L$ basis and then $E=H-2L$ (cf. also \diacrom where
the monodromy transformations correspond to
shifts in the $(H,L)$ basis not the $(E,L)$ basis)}
here is $H^6, H^4, H^0, H^2$. So the degree in the
$t$'s is $i$ on $H^{2i}$.
Here the cubic prepotential is given
for ${\bf P}^4_{1,1,2,2,2}[8]$ by\dave
\eqn\prep{{\cal F}=-{4\over 3}t_1^3-2t_1^2t_2-2t_1t_2+{7\over 3}t_1+t_2}
and for ${\bf P}^4_{1,1,2,2,6}[12]$ by
\eqn\preptw{{\cal F}=-{2\over 3}t_1^3-t_1^2t_2+{13\over 6}t_1+t_2}
with ${\cal F}^i={\partial{\cal F}\over\partial{t_i}}, i=1,2$.
It contains in the purely cubic part the properly normalized (with -1/3!)
intersection numbers $C_{ijk}$.
One then gets in a neighborhood of the large radius limit
for ${\bf P}^4_{1,1,2,2,2}[8]$ resp. ${\bf P}^4_{1,1,2,2,6}[12]$
(or of the large complex structure limit for the mirrors)
$$\Pi=\pmatrix{{2\over 3}t_1^3+t_1^2t_2+{13\over 6}t_1+t_2\cr
               {1\over 6}-2t_1t_2\cr
               1-t_1^2\cr
               1\cr
                t_1\cr
                t_2}, \ \ \
\Pi=\pmatrix{{4\over 3}t_1^3+2t_1t_2+{7\over 3}t_1+t_2\cr
               -4t_1t_2+4t_1-2t_2+{1\over 3}\cr
               -2t_1^2-2t_1+1\cr
               1\cr
               t_1\cr
               t_2}$$

\subsec{Monodromy}
The vectors just mentioned have a monodromy behaviour (represented by a
matrix $S_i$) under
$t_i \rightarrow t_i +1$ easy to read off. It belongs
to a divisor ${\cal D}_i$ in the moduli space where the corresponding
cycle becomes small. Note that
infinity is the fixpoint of the mentioned shift in the respective
variable; furthermore in the matrices $R_i:=S_i-{\bf 1}$ one finds
$R_iR_jR_k=C_{ijk}Y$ with $Y$ a matrix independent of $i,j,k$, i.e.
the $R_i$ fulfill the algebra of the $D_i$. Note also (cf. \dave\ )
that multiplication with $D_i$ corresponds
to $N_i:=log (1+R_i)$. This relation means that $ch( [D_i] )=e^{D_i}$
(where $[D_i]$ is the line bundle associated with the divisor)
corresponds to the monodromy matrix
$S_i=1+R_i$ on the Kaehler period vector.
This is the Kontsevich relation between tensoring
by the line bundle and monodromy matrix purely on the $H^{even}$ side;
mirror symmetry identifies this with the
monodromy of the type IIB periods.

Let us write down the monodromy matrices $S_i$ about the two divisors
(following the notation of \dave)
${\cal D}_H=D_{(0,-1)}$ and ${\cal D}_L=C_{\infty}=D_{(1,0)}$
for the degree 8 and the degree 12 model
$$S_L=\pmatrix{1&0&-1&2&-2&0\cr
               0&1&0&-2&-4&0\cr
               0&0&1&0&0&0\cr
               0&0&0&1&0&0\cr
               0&0&0&0&1&0\cr
               0&0&0&1&0&1}, \ \ \
S_H=\pmatrix{1&-1&-2&6&4&0\cr
               0&1&0&4&0&-4\cr
               0&0&1&-4&-4&0\cr
               0&0&0&1&0&0\cr
               0&0&0&1&1&0\cr
               0&0&0&0&0&1}$$
$$S_L=\pmatrix{1&0&-1&2&0&0\cr
               0&1&0&0&-2&0\cr
               0&0&1&0&0&0\cr
               0&0&0&1&0&0\cr
               0&0&0&0&1&0\cr
               0&0&0&1&0&1}, \ \ \
S_H=\pmatrix{1&-1&-2&5&2&1\cr
               0&1&0&0&0&-2\cr
               0&0&1&-1&-2&0\cr
               0&0&0&1&0&0\cr
               0&0&0&1&1&0\cr
               0&0&0&0&0&1}$$
The monodromy about the conifold locus is given for both models in the $(E,L)$
basis by

$$T=\pmatrix{1&0&0&0&0&0\cr
           0&1&0&0&0&0\cr
           0&0&1&0&0&0\cr
           -1&0&0&1&0&0\cr
           0&0&0&0&1&0\cr
           0&0&0&0&0&1}$$
For easier comparison let us
note that Fig. 1 of \dave\ shows among other things the divisors
(of compactified moduli space)
$C_{\infty}$ given by $\phi \ra \infty , \psi \ra \infty$,
the conifold locus $C_{con}$
given by $864\psi^6+\phi=\pm 1$ (or $(1-x)^2-x^2y=0$ in the variables
$x:=-\phi/(864\psi^6), \; y:=1/\phi^2$), $C_1$ given by $y=1$ and $C_0$
given by $\psi =0$ (with its singular point $\psi=0=\phi$); in a
heterotic language one would describe $C_{\infty}$ by $S=\infty$, $C_1$
by $S=T$ and the points of intersection of $C_{con}$ with $C_{\infty}$
resp. $C_1$ with $C_{\infty}$ by $T=i$ resp. $S=T=\infty$, leading to
the Seiberg-Witten point resp. the large radius point (or
the large complex structure point in the mirror interpretation);
'leading to' because actually and more precisely both of these points
belong to the points in the moduli space which have to be blown up to
get divisors with normal crossings only (cf. Fig. 4 of \dave\ ). In the
course of the resolution process at these two points are the divisors
$E_2$ resp. $D_{(0,-1)}$ introduced which intersect
$C_{\infty}=D_{(1,0)}$  in the point $u=\infty$ at infinity of the
Seiberg-Witten u-plane resp. in the large complex structure point.
One finds that $S_L=(ATB)^{-1}, \; S_H=(AT)^{-2}$ where $A, \; B, \;
T$ are the monodromies around $C_0, \; C_1, \; C_{con}$.
Note also that, for example,
the association between the divisor $L$, the $K3$, and
the divisor ${\cal D}_L$ in the moduli space becomes explicit
as the making the basis $P^1$ of the $K3$ fibration large is
equivalent to making the $K3$ itself hierarchically small.

Note that our $S_D$ matrices (and the $T$ matrix)
are related to those of \dave\ by
\eqn\davemat{\eqalign{m\cdot {\cal D}_D\cdot m^{-1}&= \tilde{\cal D}_D\cr
K\cdot \tilde{\cal D}_D\cdot K^{-1}&=S_D}}
with
$$m=\pmatrix{-1&1&0&0&0&0\cr
           1&3&3&2&1&0\cr
           0&1&1&1&0&0\cr
           1&0&0&1&0&0\cr
           -1&0&0&1&1&0\cr
           2&0&0&-2&1&1}, \ \ \ \
K=\pmatrix{1&0&0&0&0&0\cr
           0&1&-3&0&0&0\cr
           0&0&1&0&0&0\cr
           0&0&0&1&0&0\cr
           0&0&0&0&1&0\cr
           0&0&0&0&0&1}$$

\subsec{D6-branes charges and bundle data}
We are now interested in an explicit map between the topological invariants
of the characteristic classes of the Chan-Paton sheaf $V$ and the brane
charges.\foot{Cf. the analysis of \diacrom\ for another two-parameter
Calabi-Yau and \morri\ for the quintic.}
Therefore let us
consider the BPS charge lattice which can be
identified with the middle cohomology lattice of the mirror manifold
$H^3({\hat X},{\bf Z})$ and consider the central charge associated to the
integral vector ${\bf n}=(n_6, n_4^1, n_4^2, n_0, n_2^1, n_2^2)$ which is
\eqn\zn{Z({\bf n})=n_6\Pi_1+n_4^1\Pi_2+n_4^2\Pi_3+n_0\Pi_4+
n_2^1\Pi_5+n_2^2\Pi_6.}
One the other side, one has in the large volume limit of $X$ the lattice
of microscopic D-brane charges (which is identified with the K-theory
lattice $K(X)$). Here one considers the effective charge $Q$ of a D-brane
state\foot{Note that for $N$ coincident D6-branes wrapping $X$, the
gauge field has to satisfy the Hermitian Yang Mills equations in order to
preserve supersymmetry
$F_{ij}= 0, \; \omega^2\wedge tr F= \omega^2\wedge c_1(V)=0$.
The first equation tells us that the vacuum gauge field $A$ is holomorphic
connection on a holomorphic vector bundle $V\rightarrow X$. The second
condition is the integrability condition which guarantees a unique solution
to the Donaldson-Uhlenbeck-Yau equation $g^{i\bar{j}}F_{i\bar{j}}=0$.}
$\eta$
given by the Mukai vector
\eqn\muk{Q=
ch(\eta)\sqrt{Td(X)}\in H^0(X)\oplus H^2(X)\oplus H^4(X)\oplus H^6(X)}
with the associated central charge
\eqn\zt{Z(t)=\int_Z {t^3\over 6}Q^0-{t^2\over 2}Q^2+tQ^4-Q^6.}
where $t=t_1H+t_2L$ is the generic K\"ahler class and expansion of \muk\ leads
to
\eqn\quu{Q=(r, ch_1(V), ch_2(V)+{r\over 24}c_2(X), ch_3(V)+{1\over 24}
ch_1(V)c_2(X))}
so one has
\eqn\zent{Z(Q)={r\over 6}t^3-{1\over 2}ch_1(V)t^2+(ch_2(V)+{r\over 24}c_2(X))t
-(ch_3(V)+{1\over 24}ch_1(V)ch_2(V))}

The comparison of $Z({\bf n})$ and $Z(t)$ leads then to a map between the low
energy charges ${\bf n}$ and the topological invariants of the
K-theory class $\eta$. We find for ${\bf P}^4_{1,1,2,2,2}[8]$
\eqn\chernc{\eqalign{r(V)&=n_6\cr
                     ch_1(V)&= n_4^1E+n_4^2L\cr
                     ch_2(V)&= (4n_4^1-2n_4^2+n_2^1)h+(-2n_4^1+n_2^2)l\cr
                     ch_3(V)&= -n_0-{2\over 3}n_4^1-2n_4^2}}
and for ${\bf P}^4_{1,1,2,2,6}[12]$ we get
\eqn\cherntw{\eqalign{r(V)&=n_6\cr
                      ch_1(V)&=n_4^1E+n_4^2L\cr
                      ch_2(V)&=n_2^1h+n_2^2l\cr
                      ch_3(V)&=-n_0-{1\over 3}n_4^1-2n_4^2}}
\subsec{Kontsevich's association}
A part of Kontsevich's association can
be described just on the bundle side; it is just the question of transport
between the 'fibration' basis \decomp\ and the K\"ahler period vector basis
(cf. \dave , p.12 and \diacrom  ,p. 7;
mirror symmetry is then "only" a rephrasing of the K\"ahler period
vector in $H^{even}(X)$
as complex structure period vector in IIB for the $H^3$ of mirror);
the relations \chernc\ , \cherntw\ and (4.22) express the mentioned basis
transformations. 

We find that the monodromy transformations $S_L, S_H$
correspond to the following automorphisms $M(D)$
\eqn\shift{ [V]\rightarrow [V\otimes {\cal O}_X(L)], \ \ \ \ \
 [V]\rightarrow [V\otimes {\cal O}_X(H)]}
with the topological invariants of $V$ changed
according to for a 'twisted' sheaf
$V'=V\otimes {\cal O}_X(D)$ to
\eqn\chernclvt{\eqalign{r(V')&= r(V)\cr
                        ch_1(V')&= ch_1(V)+rD\cr
                        ch_2(V')&= ch_2(V)+ch_1(V)D+{r\over 2}D^2\cr
                        ch_3(V')&= ch_3(V)+ch_2(V)D
                                         +{1\over 2}ch_1(V)D^2+{r\over 6}D^3}}
and we see that the linear transformations acting on ${\bf n}$ corresponding
to $D=H,L$ are
\eqn\lintr{M(D)=S_D^{-1}}

Let us consider now a second type of monodromy transformation  proposed by
Kontsevich. He proposed that the monodromy $T$ about the conifold locus
of the mirror corresponds to the automorphism of the derived category
whose effect on cohomology can be denoted by
(where ${\bf 1}_X$ is the standard generator of $H^0(X,{\bf Q})$)
\eqn\smat{{\cal S}:\gamma\rightarrow \gamma-\Big(\int \gamma\wedge Td(T_X)\Big)
\cdot {\bf 1}_X}
corresponding to a change in the topological invariants of $V$
\eqn\topv{ch(V)\rightarrow ch(V)-({{ch_1(V)c_2(X)}\over 12}+ch_3(V))}
using the expression of the prepotential in the large radius limit
\eqn\prepotl{{\cal F}={1\over 6}(t\cdot J)^3-{c_2(X)\over 24}(t\cdot J)+....}
where $(t\cdot J)=\sum t_aJ^a$ (in particular we have $J_1=E+2L$ and $J_2=L$)
and
using the period vector $\Pi$ we find the expression valid for both models
\eqn\modeli{\eqalign{ch_1(V)&=n_4^1(J_1-2J_2)+n_4^2J_2\cr
                     ch_3(V)&=-(n_4^1(J_1-2J_2)+n_4^2J_2){c_2(X)\over 12}-n_0}}
leading to the universal shift
\eqn\shiuft{n_6\rightarrow n_6+n_0}
comparing this to the monodromy we find that the linear transformation
acting on ${\bf n}$ corresponds to
\eqn\lina{{\cal S}=T^{-1}}

\newsec{FM-transform as Monodromy}

Perhaps most important in view of the
other investigations in the paper is the question whether the
Fourier-Mukai transform will be related to a monodromy matrix in the
sense of Kontsevich's proposal. As the corresponding matrices are by
mirror symmetry identifiable already on the bundle side (as worked out
for some two parameter examples in the foregoing section) this comes
down, in a similar elliptic example, exxentially to
the question whether the matrix $M$ in \precisematrix\ (or
the corresponding matrix in the sLag side),
is generated by $S_H, S_L, T$. This will be made precise below.

Now recall that the transformations in the Kontsevich association
considered in the foregoing section can themselves be considered as FM
transforms. This time
we will not work on the fibre product but rather on the ordinary product.
\eqn\mats{\matrix{\;\;\; X\times {\tilde X}&\buildrel q_2 \over
\longrightarrow&{\;\;\tilde X}\cr
\scriptstyle{q_1} \biggr \downarrow & &
\scriptstyle{\pi_2}\biggr \downarrow \cr
     \;\;\; X&\buildrel {\pi_1} \over\longrightarrow&\;\;B}}
One defines the Fourier-Mukai functors $S^i$
by associating with every sheaf $V$ on $X$ the complex of sheaves
$S_{\cal E}^i(V)$ on $X$ (where $X$ and $\tilde
X$ are identified as before)
\eqn\sif{S_{\cal E}^i(V)=R^ip_{1*}(p_2^*(V)\otimes {\cal E})}
where the kernel ${\cal E}\in D(X \times X)$ is an object in the derived
category. We can also define the full FM transformation at the derived
category level
\eqn\sifd{S_{\cal E}({\cal G})=Rp_{1*}(p_2^*({\cal G})\otimes {\cal E})}
For example the twist transformation
${\cal G}\ra {\cal G}\otimes {\cal L}$ considered in \shift\ comes
from ${\cal E}= {\cal O}_{\Delta}\otimes q_2^*({\cal L})$ where $\Delta$ is the
diagonal of $X\x X$. In particular $S_{{\cal O}_\Delta}$ is the identity.
The other
operation (corresponding to the 'conifold monodromy') considered before
corresponds
to an
${\cal E}$ whose cohomology is the ideal sheaf
${\cal I}_{\Delta}$ of the diagonal of $X\x X$. The topological invariants
of this
Kontsevich FM transform can be easily obtained: The exact sequence
$$
0\ra {\cal I}_\Delta \ra{\cal O}_{X\x X}\ra {\cal O}_\Delta\ra 0
$$
leads to a triangle triangle in the derived category $D(X)$
relating the corresponding FM transforms with these
kernels:
\eqn\tri{S_{{\cal I}_\Delta}({\cal G})\to S_{{\cal O}_{X\x X}}({\cal G})\to
S_{{\cal O}_\Delta}({\cal G})\to S_{{\cal I}_\Delta}({\cal G})[1]}
and then
\eqn\chK{ch (S_{{\cal I}_\Delta}({\cal G}))=ch (S_{{\cal O}_{X\x X}}({\cal
G}))-ch
(S_{{\cal O}_\Delta}({\cal G}))=ch (S_{{\cal O}_{X\x X}}({\cal G}))-ch
({\cal G}) }
Riemann-Roch gives now that $ch_i(S_{{\cal O}_{X\x X}}({\cal G}))=0$ for
$i>0$ and
then
\eqn\chKK{ch (S_{{\cal I}_\Delta}({\cal G}))=ch_0(S_{{\cal O}_{X\x X}}({\cal
G}))-ch  ({\cal G}) =\Big(\int ch({\cal G})\wedge Td(T_X)\Big)-ch({\cal G})}
so that we recover the ``gamma shift'' \smat\ as
expected when $\cal G$ is a WIT$_1$ sheaf.

Now our fibrewise FM transform is specified by using
${\cal E}=j_*{\cal P}$ where $\cal P$ is
the Poincar\'e of \poinc\ and $j: X\times _B X \ra X \x X$ is the natural
embedding\foot{The relationship between fibrewise FM and full FM is that
(cf. also
\Thomas\ ) FM on $X\times_B X$ with respect to any sheaf $P$
is the same as FM on $X\times X$ with respect to $j_*P$ as
is clear from the fact that $p_i=q_i\circ j$ implies
$R p_{1*}=R(q_{1*}\circ j)=(R q_{1*})\circ j$
(the closed immersion $j$ has no higher derived images) and then
$R p_{1*}(p_2^*(V)\otimes P)=R q_{1*}(j_*(p_2^*(V)\otimes P))=
R q_{1*}(j_*(j^*(q_2^*(V))\otimes P))=R q_{1*}(q_2^*(V)\otimes P)$
by the projection formula.
}, so it is build up out of the structures given by the divisors
(resp. their associated line bundles)
$\sigma$ and
$\pi^*K_B$ on the one hand and by ${\cal O}_{\Delta}$ on the
other.

Since we have the $M$ matrix equation when we transform
bundles $V$ with the inverse FM, we concentrate on it. This
FM is ${\hat S}=S_{j_*{\cal Q}}$ with
$$ {\cal Q}={\cal P}^*\otimes q^*(K_B^{-1})={\cal I}
\otimes p_1^*{\cal O}(\sigma)\otimes p_2^*{\cal O}(\sigma)\otimes q^*K_B^{-2}
$$
where ${\cal I}$ is the ideal of the diagonal of $X\times_B X$. Then
\eqn\factor{{\hat S}=(\otimes q^*K_B^{-2})\circ(\otimes{\cal O}(\sigma))\circ
S_{j_*{\cal I}}\circ (\otimes {\cal O}(\sigma))}
or
$$
{\hat S}=S_{{\cal O}_\Delta(2c_1)}\circ S_{{\cal O}_\Delta(\sigma)}\circ
S_{j_*{\cal I}}\circ S_{{\cal O}_\Delta(\sigma)}
$$
We have then written our fibrewise FM ${\hat S}$ as a
composition of three Kontsewich full FM transforms and one FM transform
$S_{j_*{\cal I}}$.

Technically one can proceed now in two ways.
One can read off the relevant
matrix from our earlier treatment of the relative FM transform in
section 2. Alternatively, and presented first,
we are going to describe
$S_{j_*{\cal I}}$ in terms of the Kontsevich full FM transform
$S_{{\cal I}_\Delta}$, to make the closest contact to that quantity.

To this end, we use the exact sequence
$$
0\to {\cal J} \to {\cal I}_\Delta \to j_*{\cal I} \to 0
$$
${\cal J}$ being the ideal of the closed immersion $j: X\x_B X \ra X \x X$. We
have as in \tri\ a triangle in the derived category $D(X)$
$$
S_{{\cal J}}({\cal G})\to S_{{\cal I}_\Delta}({\cal G})\to
S_{j_*{\cal I}}({\cal G})\to S_{{\cal J}}({\cal G})[1]
$$
and then
\eqn\chide{
ch(S_{{\cal I}_\Delta}({\cal G})) = ch (S_{{\cal
I}}({\cal G}))+ ch(S_{{\cal J}}({\cal G}))
}
The term $ch(S_{{\cal I}_\Delta}({\cal G}))$ is given by \chKK\ whereas
$ch(S_{{\cal
J}}({\cal G}))$ can be computed from
$$
ch({\cal J})=1-j_*(1)+{1\over2}j_*(1)\cdot
q_1^*(c_1)-{1\over 2}
j_*(1)\cdot q_1^*(c_1^2)=1-j_*(1)+{1\over2}j_*(1)\cdot
q_2^*(c_1)-{1\over 2}
j_*(1)\cdot q_2^*(c_1^2)
$$
($j_*(1)$ is the class of $X\times_B X$ in $H^4(X\times
X)$).

If we write
$$
ch_0({\cal G})=n_{\cal G},\
ch_1({\cal G})=x_{\cal G}\sigma+S_{\cal G},\
ch_2({\cal G})=\sigma \eta_{\cal G}+a_{\cal G}F,\
ch_3({\cal G})=s_{\cal G}
$$
\chKK\ now reads
\eqn\chIDelta{\eqalign{ch_0(S_{{\cal I}_\Delta}({\cal G}))&=s_{\cal
G}-{1\over 12} x_{\cal
G}\sigma c_1^2+\sigma c_1 S_{\cal G}+{1\over 12}x_{\cal G}\sigma
c_2-n_{\cal G}\cr
&= (ch_3{\cal G}+{ch_1{\cal G}c_2(X)\over 12})-ch_0{\cal G}\cr
ch_1(S_{{\cal I}_\Delta}({\cal
G}))&=-ch_1{\cal G}\cr
 ch_2(S_{{\cal I}_\Delta}({\cal G}))&=-ch_2{\cal G}\cr
ch_3(S_{{\cal I}_\Delta}({\cal G}))&=-ch_3{\cal G}}
}
and $ch(S_{{\cal
J}}({\cal G}))$ is given by
\eqn\chJ{
\eqalign{ch_0(S_{\cal J}({\cal G}))&=s_{\cal G}-{1\over 12} x_{\cal
G}\sigma c_1^2+\sigma c_1 S_{\cal G}+{1\over 12}x_{\cal G}\sigma
c_2-x_{\cal G}\cr
&= (ch_3{\cal G}+{ch_1{\cal G}c_2(X)\over 12})-x_{\cal G}
\cr ch_1(S_{\cal J}({\cal
G}))&=-n_{\cal G} c_1-\eta_{\cal G}+{1\over 2}x_{\cal G}c_1\cr
 ch_2(S_{\cal
J}({\cal G}))&=({1\over 2}n_{\cal G}-{1\over 12}x_{\cal
G})c_1^2-c_1 S_{\cal G}+{1\over 2}\eta_{\cal G}c_1-s_{\cal G} F\cr
ch_3(S_{\cal J}({\cal G}))&=0}
}
This gives then the needed information about $ch(S_{\cal I}({\cal G}))$
which can alternatively be computed also from
\chois.
If one writes
$$
ch_0({\cal G})=n_{\cal G},\
ch_1({\cal G})=x_{\cal G}\sigma+S_{\cal G},\
ch_2({\cal G})=\sigma \eta_{\cal G}+a_{\cal G}F,\
ch_3({\cal G})=s_{\cal G}
$$
then
\eqn\chI{\eqalign{ch_0(S_{\cal I}({\cal G}))&=x_{\cal G}-n_{\cal G}
\cr ch_1(S_{\cal I}({\cal
G}))&=-x_{\cal G}\sigma- S_{\cal G}+(n_{\cal G}-{1\over 2}x_{\cal
G})c_1+\eta_{\cal G}\cr
 ch_2(S_{\cal I}({\cal G}))&=-\sigma\eta_{\cal G}-a_{\cal
G}F-({1\over 2}n_{\cal G}-{1\over 12} x_{\cal G})c_1^2+c_1S_{\cal
G}-{1\over 2}\eta_{\cal G}c_1+s_{\cal G}F\cr ch_3(S_{\cal I}({\cal
G}))&=-s_{\cal G}}}
Let us now start with a sheaf $V$ whose Chern character is
written as in \chV\ in the form $ch_0(V)=n$, $ch_1(V)=x\sigma+S$,
$ch_2(V)=\sigma \eta+aF$, $ch_3(V)=s$. By applying the composition
\factor\ we obtain
$$
ch({\hat S}(V))=ch({\cal O}(2c_1))\cdot ch({\cal O}(\sigma))\cdot ch(S_{\cal
I}({\cal G}))
$$
with ${\cal G}=V\otimes{\cal O}(\sigma)$.
The Chern character of ${\cal G}$ is given in terms of $ch(V)$ by \chernclvt.
We then have
\eqn\chG{\eqalign{n_{\cal G}&=n,\quad x_{\cal G}=x+n, \quad S_{\cal
G}=S ,\quad \eta_{\cal G}=\eta-{1\over2}nc_1+S-x c_1 \cr
a_{\cal G}&=a,\qquad s_{\cal
G}=s-\sigma c_1\eta+a+{1\over 2} x\sigma c_1^2 -{1\over 2}\sigma c_1 S+{1\over
6}n\sigma c_1^2}}
Now from \factor\ we expect the following relation for the monodromy
corresponding to the relative FM transform expressed as a product
of known Kontsevich monodromies and the monodromy around the diagonal
of the fibre product whose ideal sheaf is ${\cal I}$
\eqn\monko{S_V=S_{\sigma}\cdot S_{c_1}^{2}\cdot S_{\cal I}\cdot S_{\sigma}}
Let us make this more explicit by considering the degree 18 model.

\noindent
{\it Example}
\noindent

Consider the elliptic fibration given by
${\bf P}^4_{1,1,1,6,9}[18]$.
This model has been extensively studied in the context of mirror
symmetry \morr\ and in the context of
D-branes on elliptic Calabi-Yau \diacrom. Among the degree
18 hypersurfaces is
\eqn\dega{z_1^{18}+z_2^{18}+z_3^{18}+z_4^3+z_5^2=0}
At $z_1=z_2=z_3=0$ the ambient space has a singular line which intersects $X$
in a single point. The blow up of this line gives an exceptional
divisor $E=P^2$ in $X$. A second divisor $L$ (defined by the first order
polynomials) is given by the elliptic surface over a line in $P^2$
and together with $E$ generates $H_4(X,{\bf Z})$. The elliptic
fibration structure is induced by the linear system $|L|$ generated by
$z_1, z_2, z_3$ mapping $X$ to ${\bf P}^2$. The section of the fibration
is given by $B_2=E$. The homology class of the elliptic
fibre in $H_2(X)$ will be
denoted by $h=L^2$. Further intersection relations are given by
\eqn\ints{E\cdot L^2=1, \ \ E^2\cdot L=-3, \ \ L^3=0, \ \ E^3=9}
Working in the $E,L$ basis
(cf. \diacrom ) the generic K\"ahler class is given
by $J=t_1E+t_2L$ with $t_1,t_2$ coordinates on the K\"ahler moduli space.

In the degree 18 model again one has $\sigma =E=B=P^2$ and
$\pi^*K_B=3L$; note that $H=3L+E$ with a corresponding multiplicative
relation $S_E=S_H\cdot S_L^{-3}$
for the matrices (as the Chern character is multiplicative).

Here the relation becomes
\eqn\rel{S_V=S_E\cdot S_L^6\cdot S_{\cal I}\cdot S_E}
whereas the matrices are given by
\eqn\matri{\eqalign{S_H&=\pmatrix{1&-1&-3&10&9&3\cr
               0&1&0&0&0&0\cr
               0&0&1&0&-3&-1\cr
               0&0&0&1&0&0\cr
               0&0&0&1&1&0\cr
               0&0&0&0&0&1}, \ \ \
S_L=\pmatrix{1&0&-1&3&2&0\cr
               0&1&0&1&0&-1\cr
               0&0&1&0&-1&0\cr
               0&0&0&1&0&0\cr
               0&0&0&0&1&0\cr
               0&0&0&1&0&1}}}
$$S_E=\pmatrix{1&-1&0&1&0&0\cr
               0&1&0&-9&0&3\cr
               0&0&1&3&0&-1\cr
               0&0&0&1&0&0\cr
               0&0&0&1&1&0\cr
               0&0&0&-3&0&1}$$
Note that these matrices commute (cf. also \dave , p. 12) and that
$S_E=S_H\cdot S_L^{-3}$.

In order to obtain the matrices $S_{\cal I}$ and $S_{V}$
we have to use the comparison
(which has been performed for this model in \diacrom)
of the central charges $Z({\bf n})$ and $Z(Q)$
which gives the relation between the middle cohomology charges ${\bf n}$
and the cohomological invariants of the vector bundle ${\cal G}$
\eqn\cchv{\eqalign{ch_0({\cal G})&=n\cr
                   ch_1({\cal G})&=\alpha E +\beta L \cr
                   ch_2({\cal G})&=\gamma EL+\delta L^2\cr
                   ch_3({\cal G})&=\epsilon }}
where the coefficients are given in terms of the BPS charge vector {\bf n}
\eqn\coeff{\eqalign{n&=n_6\cr
                    \alpha&=n_4^1,\ \ \ \beta=n_4^2 \cr
                    \gamma&={3\over 2}n_4^1+n_2^2,\ \ \ \delta=
{3\over 2}n_4^2+n_2^1,\cr
                    \epsilon&=-n_0+{1\over 2}n_4^1-3n_4^2}}
Now the FM transforms ${\hat S}(.)$ and $S_{\cal I}(.)$ induce a linear
transformation on the BPS charge lattice (i.e. one compares the coefficients
of \cchv\ with the new coefficients given by the cohomological
invariants of ${\hat S}(.)$ and $S_{\cal I}(.)$).
$$S_{\cal I}=\pmatrix{1&0&3&-9&0&0\cr
           1&1&3&-9&-1&0\cr
           0&0&1&0&0&0\cr
           0&0&0&1&1&0\cr
           0&0&0&0&1&0\cr
           0&0&1&-3&0&1}, \ \ \
S_{V}=\pmatrix{0&-1&0&1&0&3\cr
           1&0&0&0&-1&0\cr
           0&0&0&0&-3&-1\cr
           0&0&0&0&-1&0\cr
           0&0&0&1&0&0\cr
           0&0&1&-3&-3&0}$$
with
\eqn\lin{S_{\cal I}^{-1}=P(S_{\cal I}(\cdot )), \ \ \
S_{V}^{-1}=P({\hat S}(\cdot ))}
denoting the linear transformations on the lattice ${\bf n}$.

Finally, using these matrices we can write the $M$ matrix as
(for $\alpha=n_4^1=0$)
\eqn\mata{M=l\cdot [S_{V}^{-1}]^t\cdot l^{-1}\cdot S_{td}}
where $l$ relates (cf. \coeff )
the period basis  with the 'fibration' basis
\decomp\
and $S_{td}$ represents the $Td(N)$ twist for this model (cf. appendix)
\eqn\lmb{l=\pmatrix{1&0&0&0&0&0\cr
           0&1&0&0&0&0\cr
           0&0&1&0&0&0\cr
           0&3/2&0&0&0&1\cr
           0&0&3/2&0&1&0\cr
           0&1/2&-3&-1&0&0}, \ \ \
S_{td}=\pmatrix{1&0&0&0&0&0\cr
           0&1&0&0&0&0\cr
           -3/2&0&1&0&0&0\cr
           0&-3/2&0&1&0&0\cr
           3/4&0&-3/2&0&1&0\cr
           0&3/4&0&-3/2&0&1}}
\newsec{Some analogues on the mirror side}
Let us note that \mata\ also shows that if we perform a linear
transformation on $H^*_{odd}(Y)$ (transforming in the ``fibre base'')
by $l$ the $M$ matrix naturally operates on the charge lattice as
$$\pmatrix{n \cr \alpha \cr \beta \cr \gamma \cr \delta \cr \epsilon}=M\cdot
\pmatrix{-\alpha\cr n\cr -\gamma \cr \beta \cr -\epsilon \cr \delta}$$
and thus \mata\ gives just the right {\it transport}\foot{Note as a caveat that
the $T^2$ which is given by the holomorphic elliptic curve used in the
fibrewise T-duality on the bundle side (an operation {\it inside}
$H^{even}(X)$)
is {\it not} contained in the $T^3$ used in the Strominger/Yau/Zaslow
T-duality employed to go {\it from} $X$ {\it to} $Y$;
only one $S^1$ of the elliptic
curve re-occurs in the $T^3$.} of the $M$ matrix
to $H^*_{odd}(Y)$ (where we have a priori only the period basis)
$${\matrix{\;\;\; H^{even}(X) &\buildrel M\over
\longrightarrow& H^{even}(X)\cr
\biggr \downarrow & &
     \biggr \downarrow \cr
     \; H^3(Y)&\buildrel {M} \over\longrightarrow&\;\; H^3(Y)}}$$
Thus the $l$ transformation
shows how to get the 'fibration' basis \decomp\
from the K\"ahler period vector
basis, resp. (after identification with the mirror side) how
the \decomp\ basis (when transported via identification to the mirror
side) is related to the complex structure period vector basis there
in IIB, making explicit the decomposition
in $H^{odd}(X^{mirror})$ corresponding to \decomp\
$${H^{odd}(X^{mirror})=H^{non-ell}\oplus H^{ell}}$$
Then the interesting question remains whether this decomposition
and 'duality' transformation on the middle {\it cohomology} of the mirror
is actually induced by a natural decomposition and map on the mirror
{\it space}. It would be interesting to unfold
that question about mirror transport
of a certain involution also in the context of involutions such as the
involution on $X$ discussed later in
sect. 5 in connection with the ${\bf Z}_2$ index theorem or
Poincar\'e duality or (fibrewise) complex structure conjugation.

\subsec{D6-brane moduli space and extended mirror conjecture}

Note that as the relevant object to consider on the bundle side
is a holomorphic cycle (so even-dimensional)
with a bundle over it correspondingly
the relevant object on the sLag side is a special Lagrangian
submanifold $C$
(three-dimensional) with a $U(1)$ bundle over it. According to
McLean's theorem the number of real (extrinsic) motions of $C$ inside $Y$ is
the same as the (intrinsic) number $b_1(C)$; these real moduli then
pair up with the same number of real moduli of the $U(1)$ bundle,
combining to $b_1(C)$ complex moduli.
\eqn\corr{h^1(End(V),X)=h^1(C)}
For a complex surface the complex dimension
of the moduli space of irreducible bundles (sheafs) is completely
determined by the Mukai vector $Q$
\eqn\muk{{\rm dim}_{\bf C}{\cal M}(Q)=Q^2+2}
which can be derived using the fact that we have a non-vanishing index
\eqn\inde{\chi(X,End(V))=\sum_{i=0}^2(-1)^i{\rm dim}H^i(X,End(V))}
Now this index becomes in the case of a Calabi-Yau three-fold
trivial. One has (cf. section 6)
by self-duality of $End(V)=V\otimes V^*$ and Serre duality
\eqn\rrs{\sum_{i=0}^3(-1)^i{\rm dim}H^1(X,End(V))=0}
\bigskip
\noindent
{\it $Z_2$ Index Theorem}

Let us consider manifolds $X$ with a group of non-trivial automorphisms
which extend to automorphisms of the vector bundle $V$ over it. For those
$X$ the Atiyah-Singer index theorem has a natural generalization \aty,\att,
the character-valued index theorem which describes how the zero modes of the
Dirac operator transform under this automorphism group.

Now the class of elliptic fibered $X$ carry such a symmetry due to
the involution $\tau$, the ``sign-flip'' in the elliptic fibers.
We will assume that at some point in the moduli space,
the $\tau$-invariant point, the symmetry can be lifted to an action on the
bundle. Then we can think of $\tau$ as decomposing
$H^i(End(V))$ into {\it even} and {\it odd} subspaces $H^i_e(End(V))$ and
$H^i_o(End(V))$. Now Serre duality involves multiplying by a holomorphic
three-form (which is odd) and thus maps $H^i_e(End(V))$ to
$H^{3-i}_o(End(V))$.
If one projects on the $\tau$-invariant part of the index problem one gets
\FMW
\eqn\tauin{-{1\over 2}\sum_{i=0}^3(-1)^iTr_{H^i(X,End(V))}\tau=
-\sum_{i=0}^3(-1)^i{\rm dim}H^i(X,End(V))_e}
which is a ``character valued index'' and can be effectively
computed by a fix point theorem \FMW. Now, using the fact that the components
of the fixed point set are of codimension two and orientable in the case of
elliptic fibered Calabi-Yau three-folds one gets for the sum
\eqn\inds{\sum_{i=0}^3(-1)^iTr_{H^i(X,End(V))}\tau=
\sum_i\int_{U_i}{{ch(End(V)_{i,e})-ch(End(V)_{i,o})}\over {1+e^{c_1(N)}}}
Td(U_i)}
where $End(V)_{i,e}$ denotes the restriction of the even resp. odd
subspaces of $End(V)$ to $U_i$ and $N_i$ is the normal bundle of $U_i$ in
$X$.
This leads \FMW\ to the dimension of the moduli space (of
$\tau$ invariant bundles; we will then assume that
$n_{odd}=h^{1,0}(C)=0$, cf. also \CDo )
\eqn\dimM{{\rm dim}_{\bf C}{\cal M}(Q)=I+2h^{1,0}(C)=r(V)-\sum_{j}\int_{U_j}
c_2(V)+2h^{1,0}(C)}
For a vector bundle as given in chapter 2 which has (with $x=0$)
\eqn\czwe{c_2(V)={L^2\over 2} -\sigma\eta-aF}
one gets (using the fact that $L^2c_1=c_1F=0$)
\eqn\iond{\sum_j\int_{U_j}c_2(V)=2L^2\sigma+\sigma\eta c_1-4a}
\bigskip
\noindent
{\it Example}

\noindent
Let us consider $X$ again being given by
${\bf P}^4_{1,1,1,6,9}[18]$ (cf. \morr ).
The fixed point set can be best described using the Weierstrass model of
$X$ given by
$x_5^2+x_4^3+fx_4+g=0$
where $f=f(x_1,x_2,x_3)$ has degree 12, and $g=g(x_1,x_2,x_3)$ has degree 18
and $(x_1,x_2,x_3)$ are coordinates on the base $B={\bf P}^2$.
The $\tau$ symmetry manifests itself in $x_5\rightarrow -x_5$ and the
two components of the fixed point set are given by $x_4=0, x_5\neq 0$ which is
isomorphic to a copy of the base (the section $E$ of $X$) on the one hand;
the other component is a triple cover of the base $B$ given by
$0=x_4^3+fx_4+g$. For this model the Chern-classes can be expressed as
\eqn\grdac{\eqalign{ch_1(V)&=n_4^2L\cr
ch_2(V)&=\big({3\over 2}n_4^2+n_2^1\big)L^2+
\big({3\over 2}n_4^1+n_2^2\big)EL\cr
ch_3(V)&=-n_0+{1\over 2}n_4^1-3n_4^2}}
(assuming here again $n_4^1=0$.)
Since $c_1(V)$ is non zero we get from the integrability condition
$\omega^2\wedge c_1(V)=0$ -guaranteeing a unique solution to the
Donaldson-Uhlenbeck-Yau equation-
that ( using $\omega=t_1H+t_2L$)
\eqn\dyu{\int_X (3t_1^2+2t_1t_2)n_4^2=0}
For the dimension of the moduli space we get therefore
\eqn\dimmo{h^1(End(V),X)=n_6-3n_2^2+6n_4^2+4n_2^1-2(n_4^2)^2}
\bigskip
\noindent
{\it Connection to FMW bundles}

\noindent
Let us now see which BPS vectors ${\bf n}$ describe the bundles
constructed by Friedman, Morgan and Witten \FMW. The bundles which are
invariant under the involution of the elliptic fiber have $c_1(V)=c_3(V)=0$
and $\eta\equiv c_1(B)$ mod 2 and $n$ is even. One has
\eqn\ctwo{c_2(V)=\eta\sigma-{(n^3-n)\over 24}c_1^2-{n\over 8}\eta(\eta-nc_1)}
therefore these bundles are described by BPS vectors
\eqn\bpsv{{\bf n}=(n_6,0,0,0,n_2^1,n_2^2)}
In order to get a dictionary between the BPS charges and the bundles
data in the FMW set-up we have to express
\ctwo\ in terms of the base $(E,S)$. Therefore setting
$\eta=ac_1(B)$ and $a$ odd we get
\eqn\esbase{c_2(V)=3aES-{{3(n^3-n)+9a(a-n)n}\over 8}S^2}
and comparing with
\eqn\csvz{c_2(V)=-n_2^2ES-n_2^1S^2}
leads then to the dictionary
\eqn\fMwv{\eqalign{n_2^2&= -3a\cr
                    n_2^1&= {{3(n^3-n)+9a(a-n)n}\over 8}\cr
                     n_6&=n}}
\newsec{Moduli space for D4-Branes and applications to
FM-transform and spectral covers}

We will be now interested in the dimension of the moduli space
of a D4-brane configuration on a divisor $D$ in $X$, so we consider
the embedding $i:D\rightarrow X$. Further let consider a vector bundle $E$
over $D$. The conditions for unbroken supersymmetry are now replaced
by the generalized Hitchin equations.
The associated K-theory class is now given by the torsion sheaf
$i_*E$ (being the extension of $E$ by zero to $X$). The Mukai vector is
then given by applying GRR for the embedding $i$
\eqn\grro{i_*(ch(E)Td(D))=ch(i_*E)Td(X)}

We will first compare
$Ext^1_X(i_*E,i_*E)$ and $Ext^1_D(E,E)$ where the first one can be
bigger by deformations (movements) of $D$ in $X$ and then for the case
of $E$ a line bundle $L$ on the spectral cover $C$ compare
$Ext^1_X(V,V)$ and $Ext^1_X(i_*L,i_*L)$ explicitly.

\bigskip
\noindent
{\it The Moduli-Space}

\noindent
The dimension of the associated moduli
space relevant here is given by the dimension of
$Ext^1_{D}(E,E)$ respectively
$Ext^1_{X}(i_*E,i_*E)$. One can in
general expect that the dimension of the moduli space associated to
$i_*E$ living over $X$
is bigger then the dimension of the moduli space of $E$ over $D$. This is
because, naively speaking, $D$ can move inside $X$ and therefore
leads to additional
deformations (the number of global deformations of $D$ in $X$) which are
related to
the number of sections of the normal bundle, i.e. the
dimension of $H^0(N)$ (cf. in the context of $F$-theory \acm).
This additional deformations play an important 
role in the comparision of D-brane moduli with the number of CFT moduli
as pointed out in the $K3$-fibration case \kllw. 
The naive picture can be made precise by considering
the long exact sequence (first written down
and proven in \thom )
\eqn\seq{0\rightarrow Ext^1_D(E,E)
          \rightarrow Ext^1_{X}(i_*E,i_*E)\rightarrow Ext^0_D(E, E\otimes N)
          \rightarrow Ext^2_D(E,E)\rightarrow}
The above exact sequence can be derived from Grothendieck
duality for the closed immersion $i: D\ra X$ (see \RD\ Section \S 6). One has
an isomorphism in the derived category
$$
R\, Hom_X(i_*E,i_*E)= R\, Hom_D(i_*E,i^!(i_*E))
$$
and $i^!(i_*E)$ is determined by the equation $i_*(i^!(i_*E))=R\,{\cal
H}om_{{\cal O}_X}(i_*{\cal O_D},i_*E)$ (where ${\cal
H}om$ stands for the Hom-sheaf). Form the exact sequence
$$
0\ra {\cal O}_X(-D)\ra{\cal O}_X\ra i_*{\cal O_D}\ra 0
$$
we read that $R{\cal
H}om_{{\cal O}_X}(i_*{\cal O_D},i_*E)$ is represented by the complex
$$
i_*E\buildrel d=0 \over \longrightarrow {\cal
H}om_{{\cal O}_X}({\cal O}_X(-D),i_*E)
$$
that is,
\eqn\rdual{i^!(i_*E) = \{\, E\buildrel d=0 \over \longrightarrow E\otimes N
\,\} }
in the derived category. Then, since $E$ is a vector bundle, we have
\eqn\extdual{\eqalign{R\, Hom_X(i_*E,i_*E)&=R\, Hom_D(i_*E, \{ \,E\buildrel
d=0
\over
\longrightarrow E\otimes N\,\} )\cr
&=R\,\Gamma(D, \{\, End(E)\buildrel d=0 \over
\longrightarrow End(E)\otimes N\,\})
}}
Taking into account the natural isomorphisms $Ext_D^i(E,E)=H^i(D,End(E))$
and $Ext_D^i(E,E\otimes N)=H^i(D,End(E)\otimes N)$, the exact sequence
\seq\ is identified with the sequence of the low terms
$$
0\ra E_2^{1,0}\ra H^1(M)\ra E_2^{0,1}\buildrel{d_2}\over\longrightarrow
E_2^{2,0}\ra H^2(M)
$$
of the spectral
sequence approaching $Ext_X^i(i_*E,i_*E)$ from the double complex associated
to an injective resolution of $\{\, End(E)\buildrel d=0 \over
\longrightarrow End(E)\otimes N\,\}$. We have here a more complete
information: since
$d=0$ the first differential of this double complex is zero and then
$d_2: E_2^{0,1}\ra E_2^{2,0}$ is zero as well. The sequence \seq\ takes now
the form
\eqn\seqtwo{0\rightarrow Ext_D^1(E,E)=H^1(D,End(E))
          \rightarrow Ext^1_{X}(i_*E,i_*E)\rightarrow H^0(D,End(E)\otimes
N) \rightarrow 0}
When $E=L$ is a line bundle on a spectral cover $D=C$, $End(L)={\cal
O}_C$ and we have
$$
0\rightarrow H^1(C,{\cal O}_C)
          \rightarrow Ext^1_{X}(i_*L,i_*L)\rightarrow H^0(C,{\cal O}_C)
\rightarrow 0
$$
so that
\eqn\dimext{\eqalign{ {\rm dim} \ Ext^1_D(L,L)&=h^1(C,{\cal O}_D)
=h^{(0,1)}(C)\cr
{\rm dim} \ Ext^1_{X}(i_*L,i_*L)&=h^1(C,{\cal O}_D)
+h^0(C,N)=h^{(0,1)}(C)+h^{(2,0)}(C)}}
where the very last formula uses that the ambient $X$ is CY.

If we consider the vector bundle $V=S^0(i_*L)$ derived from $L$ by
the spectral cover construction (or by the FM transform), then $V$ is WIT$_1$
and its unique inverse FM transform goes back to
${\hat S}^1(V)=i_*L$.
Then, by the ``Parceval isomorphism'' (see \Muk\ , \Mac\ or \RPo\ ), we
have
\eqn\parc{Ext^1_X(V,V)=Ext^1_X({\hat S}^1(V),{\hat S}^1(V))
=Ext^1_X(i_*L,i_*L)}
and  then
\eqn\dimextv{\eqalign{ {\rm dim} \ Ext^1_D(L,L)&=h^1(C,{\cal
O}_D)=h^{(0,1)}(C)\cr {\rm dim} \ Ext^1_{X}(V,V)&=h^1(C,{\cal O}_D)
+h^0(C,N)=h^{(0,1)}(C)+h^{(2,0)}(C)}}
We now want to show {\it explicitely}
that ${\rm dim} Ext^1_X(V,V)={\rm dim} Ext^1_X(i_*L,i_*L)$ without using
\parc\ : Let us recall that from Serre duality one has
\eqn\see{\sum_{i=0}^3(-1)^i {\rm dim}\; Ext^i(V,V)=0}
for $V$ on a Calabi-Yau. However, since we work on elliptic
CY we can use the character valued index and compute in the spectral
cover representation that (cf. \CDo\ )
\eqn\cvi{{\rm dim}\; H^1(X,End(V))=h^{(2,0)}(C)+h^{(1,0)}(C)}
showing together with \dimext\ the wished for agreement.

\bigskip
\noindent
{\bf Acknowledgements}: We would like to thank G. Hein, R. Thomas and
C. Vafa for discussions. The third author is a member
of VBAC
 (Vector Bundles on Algebraic Curves),
which is partially supported by EAGER (EC FP5 Contract no.
HPRN-CT-2000-00099). He also acknowledge support under grant BFM2000-1315
 of the spanish DGI.

\centerline{\bf Appendix}
\noindent
In this appendix we illustrate why in the three-fold case
the sole use of the T-functor known from the $K3$-case
to map the Chern classes of the bundle and its dual
is insufficient to exhibit as transformation matrix
the adiabatic extension \precisematrix\
of the usual T-duality matrix \ellmatrix\
on the fibre.

Let us recall the findings in this case (we phrase them here in the
language of \FMW\ ).
The Chern character of $V$ is given by\foot{Note that $ch_2(V)=-c_2(V)$
can be computed \FMW\ from its restriction to $\sigma=B$ via
$(\pi|_C)_*(e^{c_1(L)}Td(C))=ch(V|_{\sigma})Td(B)$ where then the term
$-\eta c_1$ is corrected/lifted to $\eta \sigma$ leading to the correction
$-(\eta \sigma +\eta c_1)$ above. For the different evaluations of
$ch_3(V)=c_3(V)/2$ cf. \cur ; note also that
$C\sim n\sigma +\eta$, $c_1(C)=-C|_C, c_2(C)=c_2(X)|_C+C^2,
c_1(L)={C+c_1 \over 2}+\gamma$ and $c_2(X)=12\sigma c_1+c_2+11c_1^2$.}
$${\eqalign{r(V)&= n\cr
                   ch_1(V)&= 0\cr
                   ch_2(V)&= -\eta \sigma -\eta c_1+
                         \pi_*({{c_1(C)^2+c_2(C)}\over 12}+{{c_1(L)
                          c_1(C)}\over 2}+{{c_1(L)^2}\over 2})-
                          {{n(c_1^2+c_2)}\over 12}\cr
                     &= -\eta \sigma +
                         \pi_*({{C^2}\over 24}+{1\over 2}\gamma^2)
                         -{n\over 24}c_1^2\cr
                   ch_3(V)&=\lambda\eta(\eta-nc_1)\cdot \sigma
                           =-\gamma C\cdot \sigma }}$$
Using the decomposition of the cohomology we find
$$
ch(V)=\pmatrix{n\cr
                 0\cr
                 0\cr
                 -\eta\cr
                 \pi_*({C^2\over 24}+{1\over 2}\gamma^2)-{n\over 24}c_1^2\cr
                 -\gamma C}
$$
Using GRR we get\foot{Note that the intersection products
in $C$ pushed down to $B$ by $\pi_*$ (as used in the computations for $V$)
equal the same intersection products pushed forward to $X$ via $i_*$ (as used
in the computations for $V~$) as only
those triple intersections survive which have one $\sigma$ factor and then the
remaining intersection product of two classes in $B$ equals the expression
obtained on the first route
(also note that $i_*1=C$, $\pi_*\sigma=\eta-nc_1$ and $\pi_*1=n$).}
on the other hand
for $ch(i_*L)$
$${\eqalign{r(i_*L)&=0\cr
                  ch_1(i_*L)&=C\cr
                  ch_2(i_*L)&=i_*(c_1(L)+{c_1(C)\over 2})\cr
                            &=i_*({c_1\over 2}+\gamma )\cr
                  ch_3(i_*L)&=i_*({C^2+3c_1^2\over 24}
                              +{\gamma (c_1+\gamma)\over 2})}}$$
\bigskip
\noindent
{\it The T-functor, analogous to $K3$, in the three-fold case}

If one tries for the T-functor again
$${T^?(\cdot)=S(\cdot)\otimes \pi^* K_B^{-1/2}}$$
one will get not completely the vector needed for the M matrix:

The Chern-characters of $T^?(V)=\tilde{V}\otimes K_B^{-1/2}
=i_* L \otimes K_B^{+1/2}$ are by \chernclvt\
given by
$${\eqalign{r(T^?(V))&=0\cr
                   ch_1(T^?(V))&=C=n\sigma+\eta\cr
                   ch_2(T^?(V))&=\gamma C\cr
                   ch_3(T^?(V))&=i_*({C^2\over 24}+{\gamma^2\over 2})}}$$
and so
$$ch(T^?(V)=\pmatrix{0\cr
                  n\cr
                  \eta\cr
                  0\cr
                  \gamma C\cr
i_*({C^2\over 24}+{\gamma^2\over 2})}$$
showing the mismatch of $nc_1^2/24$ between the
fifth entry of $ch(V)$ and the sixth entry of
$ch(T^?(V))$. Even without explicit computation of both sides one can
see directly why this can not work as
(the last step uses $V|_B={\pi_C}_*L$
following from $V={\pi_1}_*(\pi_2^*L\otimes {\cal P})$)
$${\eqalign{
\int_X ch(i_*L)Td(X)&=\int_X i_*(ch(L)Td(C))=\int_C ch(L)Td(C)\cr
                    &=\int_B {\pi_C}_*(ch(L)Td(C))
                     =\int_B ch({\pi_C}_*L)Td(B)\cr
                    &=\int_B ch(V|_B)Td(B)}}$$
This shows that
$${\int_X ch_3(i_*L) + C{c_2(X)\over 12}=
\int_Bch_2(V|_B) +n {c_1^2+c_2\over 12}}$$
or, in other words\foot{here $|_{(2)}$ denotes the terms of
complex dimesion 2; note that it is $-(Td(X)-Td(B))|_{(2)}$
what is used and not $-{Td(X)\over Td(B)}|_{(2)}$;
the latter would have a ${c_1^2\over 4}$ correction.}
(where in the last step $-\omega_{ch_2(V)}=ch_2(V|_B)-\eta c_1$ is
used\foot{with the notation
$c_2(V)=\eta\sigma + \omega$ where $\omega\in H^4(B)$ (pullback for
$\omega$ understood)})
$${\eqalign{\int_X ch_3(T^?(V))
&= \int_B ch_2(V|_B)+\Bigl (ch K_B^{1/2}-Td(X)+Td(B)\Bigr ) |_{(2)} \cdot C\cr
&= \int_B ch_2(V|_B)+ (-{c_1^2\over 8}-{c_2(X)\over
12}+{c_1^2+c_2\over 12}) \cdot C\cr
&= \int_B ch_2(V|_B)+ n{c_1^2\over 24}-\eta c_1)\cr
&= -\omega_{ch_2(V)}+ n{c_1^2\over 24}}}$$
\bigskip
\noindent
{\it The amended T-functor}

However, if we introduce the T functor
$${T(\cdot)=S(\cdot )\times Y}$$
(so that for example
$T(V)=\tilde{V}\otimes Y$ and $T(i_* L)=V\otimes Y$)
with a sheaf $Y$ of
$${ch(Y)=(1+{c_1\over 2}+{c_1^2\over 6})
=1/(1-{c_1\over 2}+{c_1^2\over 12})=Td(N)^{-1}}$$
where $N$ denotes the normal bundle of $B$ in $X$, then we find
$$
ch(T(i_* L))=ch(V)ch(Y)=\pmatrix{0&1&0&0&0&0\cr
-1&0&0&0&0&0\cr
0&0&0&1&0&0\cr
0&0&-1&0&0&0\cr
0&0&0&0&0&1\cr
0&0&0&0&-1&0}ch(i_*L)
$$
Equivalently we can write
(note that $Td(N)=ch(K_B^{1/2})\cdot (1-{c_1^2\over 24})$)
$${\eqalign{ch(V)&=M\cdot ch(i_*L)Td(N)\cr
                       &=M\cdot ch(i_*L)(1-{c_1\over 2}+{c_1^2\over 12})
=M\cdot ch(i_*L)(1-{c_1\over 2}+{c_1^2\over 8}-{c_1^2\over 24})}}$$
But note that the appearance of the $Td(N)^{-1}$
term can not be understood by
taking for example just $Y=j_*({\cal O}_{\sigma})$, as
GRR for the embedding $j:B\hookrightarrow X$ gives actually
$${\eqalign{
ch(j_*({\cal O}_{\sigma}))&
=j_*\big( ch(({\cal O}_{\sigma}))Td(B)\big) Td(X)^{-1}
=j_*\big( ch(({\cal O}_{\sigma}))Td(B)/j^*Td(X)\big) \cr
          &=j_*\big( ch(({\cal O}_{\sigma}))Td(N)^{-1}\big)
           =j_*(Td(N)^{-1})=\sigma Td(N)^{-1}}}$$

\noindent
{\it Td(N) twist as matrix}
\noindent

One can present the $Td(N)$ twist
as a matrix if one considers its operation on a general cohomology vector
$$
v=\pmatrix{n\cr
             x\cr
             S\cr
             \eta\cr
             a\cr
             s}=n+(x\sigma +S)+(\eta\sigma + aF)+s
$$
that for a twist by $Td(N)=1-{c_1\over 2} +{c_1^2\over 12}$ (or similar for
$Td(N)^{-1}=1+{c_1\over 2} +{c_1^2\over 6}$)
$$
{Td(N)v
=n+(x\sigma +S)+(\eta\sigma + aF)+s
=\pmatrix{1&0&0&0&0&0\cr
          0&1&0&0&0&0\cr
-{c_1\over 2}&0&1&0&0&0\cr
0&-{c_1\over 2}&0&1&0&0\cr
{c_1^2\over 12}&0&-{c_1\over 2}&0&1&0\cr
0&{c_1^2\over 12}&0&-{c_1\over 2}&0&1}v}
$$
which when one also decomposes the base
cohomology becomes a matrix of numbers.

\listrefs
\bye